\documentclass[11pt,a4paper]{amsart}
\usepackage[english]{babel}
\usepackage[utf8]{inputenc}
\usepackage[foot]{amsaddr}
\usepackage{microtype}
\usepackage{amsmath}
\usepackage{amsfonts}
\usepackage{bbm}
\usepackage{amssymb,amsthm}
\usepackage[format=plain, indention=.5cm]{caption}
\usepackage{subcaption}
\usepackage{graphicx}
\usepackage{color}
\usepackage[usenames,dvipsnames]{xcolor}
\usepackage{textcomp}
\usepackage[a4paper]{geometry} 
\usepackage[english]{babel}
\usepackage{mathtools}
\usepackage[
final,
colorlinks,
linkcolor = BrickRed,
citecolor = OliveGreen,
filecolor = black,
urlcolor = BrickRed,
bookmarks=true,
bookmarksnumbered=true,
bookmarksopen=true,
bookmarksopenlevel = 0,
hypertexnames=false,
]{hyperref}
\usepackage[
style=alphabetic,
giveninits=true,
maxbibnames=6,
minbibnames=6,
]{biblatex}

\renewbibmacro*{doi+eprint+url}{%
	\printfield{doi}%
	\newunit\newblock
	\iffieldundef{doi}
	{\usebibmacro{url+urldate}}
	{}%
	\newunit\newblock
	\usebibmacro{eprint}%
}

\usepackage[capitalise,noabbrev,nameinlink]{cleveref}
\usepackage{siunitx}
\newtheorem{thm}{Theorem}[section]

\newtheorem{prop}[thm]{Proposition}

\newtheorem{rem}[thm]{Remark}
\graphicspath{{Pictures/}}
\usepackage{float}
\usepackage{multirow}
\usepackage{xcolor}
\usepackage{tikz}
\usetikzlibrary{patterns}
\usepackage{pgfplots}
\usepackage{booktabs}
\usepackage{comment}

\newcommand{\N}{\mathbb{N}}
\newcommand{\R}{\mathbb{R}}
\newcommand{\Z}{\mathbb{Z}}
\newcommand{\diag}{\operatorname{diag}}
\newcommand{\tr}{\operatorname{tr}}
\newcommand{\adj}{\operatorname{adj}}

\newcommand{\id}[1][{}]{\operatorname{id}_{#1}}

\newcommand{\sign}{\operatorname{sign}}
\newcommand{\cc}{{\operatorname{c}}}
\newcommand{\pc}{{\operatorname{pc}}}

\newcommand{\Pid}{\Pi_d}
\newcommand{\Loss}{\mathcal{L}}

\newcommand{\SOd}{\mathcal{SO}(d)}
\newcommand{\M}{\mathcal{M}}
\newcommand{\m}{m}
\newcommand{\dt}[1][t]{\,\mathrm{d}#1}
\newcommand{\dx}{\dt[x]}

\newcommand{\pred}{{\rm pred}}
\newcommand{\sym}{{\rm sym}}
\newcommand{\ineq}{{\rm ineq}}
\newcommand{\mse}{{\rm mse}}
\newcommand{\shift}{{\rm shift}}
\definecolor{unia-green}{RGB}{0, 101, 97}
\definecolor{unia-pink}{RGB}{173, 0, 124}
\definecolor{unia-yellow}{RGB}{246, 168, 0}
\definecolor{unia-orange}{RGB}{235, 105, 11}
\definecolor{unia-red}{RGB}{212, 0, 45}
\definecolor{unia-lightblue}{RGB}{0, 174, 207}
\definecolor{unia-blue}{RGB}{0, 135, 193}
\definecolor{unia-lightgreen}{RGB}{72, 147, 36}
\colorlet{coltrainLoss}{unia-green}
\colorlet{colvalLoss}{unia-orange}
\colorlet{colvalLossMse}{unia-lightgreen}
\colorlet{colvalLossSym}{unia-red}
\colorlet{colvalLossIneq}{unia-blue}
\colorlet{colPhi}{gray}
\colorlet{colPhiPc}{unia-lightgreen}
\colorlet{colPhiPcSVPCLP}{unia-lightgreen}
\colorlet{colPhiPcPred}{unia-pink}

\title[Neural Network Polyconvexification]{Neural Network Enhanced Polyconvexification of Isotropic Energy Densities in Computational Mechanics}
\author[]{L.~Balazi$^{*}$, T.~Neumeier$^{*}$, M.~A.~Peter$^{\dagger}$, D.~Peterseim$^{\dagger}$}
\address{${}^{*}$ Institute of Mathematics, University of Augsburg, Universit\"atsstr.~12a, 86159 Augsburg, Germany}
\address{${}^{\dagger}$ Institute of Mathematics \& Centre for Advanced Analytics and Predictive Sciences (CAAPS), University of Augsburg, Universit\"atsstr.~12a, 86159 Augsburg, Germany}
\email{\{loic.balazi, timo.neumeier, malte.peter, daniel.peterseim\}@uni-a.de}
\date{\today}

\begin{document}

\thanks{The authors gratefully acknowledge funding from the German Research Foundation (DFG) within the Priority Programme 2256 \emph{Variational Methods for Predicting Complex Phenomena in Engineering Structures and Materials} (project number 441154176, reference IDs PE1464/7-2 and PE2143/5-2). Furthermore, we would like to thank the Bavarian State Ministry of Science and the Arts for funding the Augsburg AI Production Network as part of the High-Tech Agenda Plus.}

\begin{abstract}
    We present a neural network approach for fast evaluation of parameter-dependent polyconvex envelopes, which are crucial in computational mechanics. Our method uses a neural network architecture that inherently encodes polyconvexity in the main variable by combining a feature extraction layer that computes the minors function on the signed singular value characterisation of isotropic energy densities with a Partially Input Convex Neural Network (PICNN). 
    The envelope inequality is weakly enforced by penalisation during training, as are the symmetries of the function. 
    As a guiding example, we focus on a pseudo time incremental variational damage problem, which is parameter-dependent on previous time-step iterates, the deformation gradient and the internal variable. 
    This problem is reformulated in terms of signed singular values and a splitting approach is applied to reduce the dimension of the parameter space, thereby making training more tractable.
    Numerical experiments show that the networks achieve favourable accuracy for engineering applications while providing high compression and significant speed-up over traditional polyconvexification schemes.     
    Most importantly, the network adapts to varying physical or material parameters, enabling real-time polyconvexification in large-scale computational mechanics scenarios.
\end{abstract}
    
\maketitle

{\tiny {\bf Key words.} Polyconvexity, Input Convex Neural Network, relaxation, isotropic damage model, parameter-dependent}\\
\indent
{\tiny {\bf AMS subject classifications.} {\bf 49J45}, {\bf 49J10}, {\bf 74G65}, {\bf 74B20}, {\bf 68T07}, {\bf 74A45} 
} 

\section{Introduction}
Many relevant problems in solid mechanics seek to identify the minimising deformations \(u\colon \Omega \to \R^d\) of a body \(\Omega\) in spatial dimension \(d \in \{2,3\}\) with respect to the energy functional 
\begin{equation} \label{eq:energyfunctional} 
I(u) = \int_{\Omega} W(\nabla u; \zeta) \dx, 
\end{equation} where the energy density \(W \colon \R^{d \times d} \times \R^{p} \to \R_{\infty} \coloneqq \R \cup \{\infty\}\), with $p\in \N$, is defined over a suitable class of admissible functions. The vector \(\zeta \in \R^p\) represents spatially varying parameters, including material history, damage evolution or phase information. In many established engineering models, the density \(W\) is either inherently nonconvex or develops nonconvexity as its parameters evolve, as in damage or phase transition models. Such nonconvexity not only poses significant mathematical challenges, including the potential non-existence of minimisers, but also leads to serious problems in numerical simulation, such as mesh dependence and reduced robustness.

A common approach to overcome these issues is relaxation via the (semi-)convexification of the energy density, where the term semi-convexity covers the notion of polyconvexity, quasiconvexity and rank-one convexity. By replacing the nonconvex energy density \(W\) with its (semi-)convex envelope, a relaxed functional is obtained, which guarantees the existence of minimisers and allows for the identification of bounds of infimising sequences for the original nonconvex problem. On the computational side, these favourable properties  translate into robust and mesh-insensitive numerical results. 

As shown in \cite{Bal76,Bal77,Bal02}, the notion of polyconvexity is well suited for nonlinear elasticity. In general, for a given function $W$, it is not straightforward to verify whether $W$ is polyconvex, which has led to extensive investigations of sufficient and necessary conditions in the literature, see e.g.\ \cite{RosSim95, Ros98, Sil99, HarNef2003, DacMar2006, CarDacHanLab1988, ButDacGan94, Kri2000, BouKreSch19, Steigmann2003, Mie05, WiePet2026, WolHolNef2026}.
However, for the relaxation of nonpolyconvex functions, explicit analytical representations of the polyconvex envelope \(W^{\pc}\) are typically not available, necessitating computational relaxation methods.
Conventional approaches for relaxation include optimisation-based and computational geometry methods, such as algorithms approximating the rank-one convex envelope \cite{Dol99, DolWal00, Bar04, BalKohNeuPetPet23, KNPPB24, KNMPPB2022, ConDol18} or dedicated algorithms for computing the polyconvex envelope \cite{Bar05,EneBosGri13,BEG15,FanHenKruMurWei2026}. However, these methods are often subject to the curse of high-dimensionality, since they require the mesh-based discretisation of the (\(d \times d\))-dimensional space of deformation gradients.
Polyconvexification algorithms further increase the dimensionality issue because they perform the convexification in the minors space of the deformation gradient. For isotropic functions, this dimension can be reduced by characterising the polyconvexity via signed singular values, as in \cite{WiePet2026} (see \cite{NPPW2024} for a corresponding polyconvexification algorithm). However, despite their efficiency, these algorithms remain resource and time consuming, especially in spatial dimension \(d=3\). As a result, conventional relaxation algorithms are often impractical for engineering applications that require iterative computations and multi-query evaluations of the polyconvex envelope for varying parameter values.

For example, frameworks such as the isotropic pseudo-time incremental damage model introduced in \cite{BalOrt2012} involve a (\(d\times d + 1\))-parameter-dependent family of convexification problems, consisting of the deformation gradient and the scalar internal variable from the previous time step that captures the material history (see \cref{sec:damage}). Such parameter dependence introduces additional dimensionality as the material parameter \(\zeta \in \R^{p}\) (for the damage example it holds \(p = d^2 + 1\)) can vary at any material point in \(\Omega\). This requires the polyconvex envelope to be computed from scratch for each parameter configuration, resulting in a (\(d\times d+ p\))-dimensional problem. Overall, the accurate approximation of the polyconvex envelope for parameter-dependent energy densities is of great importance for the practical application of computational relaxation techniques.

To overcome the high computational cost of conventional relaxation schemes and make them practically feasible, we propose a neural network-based approach that compresses the parameter-dependent polyconvex envelope into the parameters of a small-scale artificial neural network that incorporates intrinsic properties of the polyconvexification problem into its architecture. Our method exploits a design that encodes polyconvexity by combining a feature extraction layer which computes the minors function based on a signed singular value characterisation of isotropic functions, with Input Convex Neural Networks (ICNNs) \cite{ICCN}. In addition, symmetry properties and envelope inequality are systematically enforced during training through penalty terms in the loss functional, ensuring that the network satisfies the necessary physical and mathematical constraints. This tailored approach not only speeds up the approximation of the polyconvex envelope, but also increases the reliability of the approximations in practical engineering applications.

Recent advances in neural networks, particularly through the development of ICNNs \cite{ICCN}, have enabled architectures that enforce convexity in energy density functions, a key requirement for accurate material modelling. For example, \cite{Faisal, Faisal1} use ANNs to learn constitutive laws from stress--strain data, while preserving fundamental principles of solid mechanics through ICNNs that ensure material stability. It should be emphasised that in these two works the ANNs are convex with respect to the strain tensor \(E=\frac{1}{2}(F^TF - \mathbb{I})\) instead of the deformation gradient $F$. These models are thus not guaranteed to be polyconvex and are therefore not elliptic for all states. In addition, \cite{klein2022} and \cite{Klein} propose machine-learning-based constitutive models for finite deformation and electro-mechanically coupled behaviour respectively, using ICNNs to enforce hyperelasticity, anisotropy and polyconvexity. In \cite{Linden}, a neural network-based hyperelastic constitutive model is introduced that inherently satisfies standard constitutive conditions using ICNNs, while \cite{Vijayakumaran} designs multi-scale heterogeneous structures with spatially varying microstructures by incorporating physical principles such as polyconvexity, objectivity and thermodynamic consistency. 

Without using ICNNs, some works have also explored ANNs for modelling the behaviour of mechanically sound materials \cite{Linka, Fernandez, Vahidullah, Vahidullah2, chen2022, FerFriWee2022}. More closely to our approach, other works focus on (Partially) Input Convex Neural Network formulations in the invariants of the right Cauchy--Green tensor \(C\), sufficient but not necessary for polyconvexity, see e.g.~\cite{KleRotValWee2023, ZheKocKum2024, KleHosKikKanRudGil2025}. Recently, \cite{GeuKurWieMos2025} introduced convex signed singular value neural networks (CSSV-NN) for isotropic hyperelastic energy formulations and in \cite{GeuKurWieZlaCanMos2026} attention was drawn to the representation of incompressible materials. 
Unlike \cite{GeuKurWieMos2025}, which aims to learn polyconvex hulls directly from values of potentially nonpolyconvex functions, our approach learns and compresses polyconvex envelopes from reference solutions obtained either analytically or via established numerical algorithms \cite{NPPW2024}. 
This ensures polyconvex data as input and provides stability in the training process, as well as accuracy and control over the learning target within the prescribed domain.
The approach aims to construct a compressed representation of the polyconvex envelope on the entire signed singular value space, facilitating its efficient evaluation in complex engineering loading scenarios.
In addition, the method extends to the more general case of parameter-dependent classes of reference envelopes: the trained network can thus be applied to entire families of polyconvexification tasks, for example within damage simulations, rather than being limited to a single envelope.

The remaining parts of this paper are structured as follows. In \cref{sec:damage}, we introduce a guiding engineering application by presenting a parameter-dependent family of energy densities arising from an isotropic damage problem. In \cref{sec:theorypoly}, we recall the signed singular value formulation of the polyconvex envelope for isotropic functions, while \cref{sec:theoryNN} introduces the properties-preserving neural networks that enforce polyconvexity, symmetries and the envelope inequality. Finally, \cref{sec:numericalexperimentsmath,sec:numericalexperimentsisodamage} provide numerical examples in two and three spatial dimensions including both mathematical benchmarks and the engineering isotropic damage formulation, rewritten in the signed singular value framework and parameter-reduced via a splitting approach.

\section{Isotropic Damage Problem: A Guiding Example} \label{sec:damage}
The numerical simulation of the isotropic damage model of \cite{BalOrt2012} is a prime representative of a problem where parameter-dependent polyconvexification is encountered in an engineering application and serves as a guiding example for our work. 

Let $F$ denote the deformation gradient, i.e.~\(F = \nabla u\) for a sufficiently smooth deformation \(u\colon \Omega \to \R^d\) of a body \(\Omega\) in spatial dimension \(d \in \{2,3\}\) as in \eqref{eq:energyfunctional}. 
In what follows, we consider two different effective hyperelastic strain energy density functions $\psi^{0}$ in the undamaged state: the Saint~Venant--Kirchhoff (STVK) model
\begin{equation} \label{eq:STVK}
    \psi^{0}_{\rm STVK}(F) = \frac{\mu}{4} \, \lvert F^T F - \mathbb{I}\rvert^2 + \frac{\lambda}{8} \, \left(\lvert F \rvert^2 - d\right)^2,
\end{equation}
representing the special case of a linear stress--strain relation in a Lagrangian formulation, and a compressible neo-Hookean (NH) model 
\begin{equation} \label{eq:NH}
	\psi^0_{\rm NH}(F) = \frac{\mu}{2} (\tr(F^T F) - d)  - \mu \ln(\det F) + \frac{\lambda}{2} \ln(\det F)^2 ,
\end{equation}
showing a nonlinear response, which also incorporates the determinant constraint for the penalisation of negative determinants of the deformation gradient.
Here, $\lambda$ and $\mu$ denote the Lamé constants. 

The pseudo-time incremental energy density function capturing isotropic damage for the time step \(k\to k+1\) is formulated in \cite[Equation~(20)]{BalOrt2012} as 
\begin{equation} \label{eq:Wdamage}
	\begin{aligned}
		W(F_{k+1}; F_{k}, \alpha_{k}) & =  \psi(F_{k+1},p(F_{k + 1}; \alpha_{k})) - \psi(F_k,\alpha_k) \\
		&\qquad + p(F_{k + 1}; \alpha_{k})\, D(p(F_{k + 1}; \alpha_{k})) - \alpha_k \, D(\alpha_k) 
		- \overline{D}(p(F_{k + 1}; \alpha_{k})) + \overline{D}(\alpha_k).
	\end{aligned}
\end{equation}
In \eqref{eq:Wdamage}, the scalar \(\alpha_k \in \R\) and the matrix \(F_k \in \R^{d\times d}\) denote the internal variable and the deformation gradient from the previous pseudo-time step \(k\), respectively, and serve as parameters for the current step \(k+1\); they are considered fixed and known. The strain energy density function \(\psi\) is given by 
\begin{equation} \label{eq:psi1minusD}
	\psi(F, \alpha) = (1 - D(\alpha)) \, \psi^{0}(F),
\end{equation}	
where \(\psi^{0}\) is an isotropic energy density such as the STVK \eqref{eq:STVK} or NH \eqref{eq:NH} examples above. Through this formulation, the characteristics of the undamaged energy density, such as isotropy or determinant constraints, enter the model. The function \(D\colon \R \to [0,1)\) is a non-decreasing exponential damage function, as considered in \cite{Mie:1995:dcd}, which maps the internal variable \(\alpha\) to the interval \([0,1)\), with \(0\) corresponding to the undamaged state and values close to \(1\) indicating full damage. Specifically, 
\begin{equation} \label{eq:damagefunction} 
	D(\alpha) = d_\infty \left(1 - \exp\left(- \frac{\alpha}{d_0}\right)\right), 
\end{equation} 
where \(d_\infty \in (0,1)\) is the \emph{asymptotic limit} (the maximum possible damage) and \(d_0 \in \R^+\) is the \emph{damage saturation} parameter. The reciprocal of \(d_0\) reflects the rate at which the asymptotic limit \(d_\infty\) is reached. The function 
\begin{equation*} \overline{D}(\alpha) = d_\infty \left(\alpha + d_0 \, \exp\left(-\frac{\alpha}{d_0}\right)\right) 
\end{equation*}
is the antiderivative of \(D\). The evolution of the internal variable \(\alpha_{k+1}\), and hence the damage \(D(\alpha_{k+1})\), is determined by a \emph{path function} \(p\), defined as \begin{equation} \label{eq:pathfunctionW}
	\alpha_{k+1} = p(F_{k+1}; \alpha_{k}) = \begin{cases}
		\psi^{0}(F_{k+1}) & \text{if } \, \psi^{0}(F_{k+1}) > \alpha_{k}, \\ \alpha_{k} & \text{else}.
	\end{cases}
\end{equation}
The function \(W\) in \eqref{eq:Wdamage} is nonconvex because the damage evolution process is governed by the path function \(p\). 
Theoretically, nonconvexity prevents the existence of minimisers. Numerically, it poses challenges such as stability problems, nonconvergence and oscillatory behaviour in the simulation of boundary value problems; see \cite{BalOrt2012,BalKohNeuPetPet23,KNPPB24,KNMPPB2022}, where these issues have been addressed by rank-one convexification, for further discussion of these issues.

\section{Polyconvexification of Isotropic Functions} \label{sec:theorypoly}
In this section, we present a polyconvexification approach based on a reformulation in terms of the signed singular values of the deformation gradient, which underpins our neural network design. For ease of notation, we omit the explicit parameter dependence and consider functions \(W\colon\R^{d \times d} \to \R_{\infty}\). The parameter-dependent formulation is obtained by direct analogy, by including the parameter vector \(\zeta\) and considering the function \(W(\,\cdot\,; \zeta)\).

Assuming that \(W\) is isotropic, we use its characterisation by signed singular values to reduce the dimension of the problem. In this formulation, polyconvexification of an isotropic function reduces to convexification of a function defined on the manifold determined by the minors of the signed singular values. This reformulation is advantageous because it reduces the domain from \(\R^{d \times d}\) to \(\R^d\), thereby reducing the computational cost of a grid-based representation of the polyconvex envelope. This reduced approach has been successfully applied in an efficient conventional polyconvexification algorithm \cite{NPPW2024} (see also \cref{sec:SVPC}). In the context of this work, it serves as the basis for learning and predicting the polyconvex envelope via properties-preserving artificial neural networks.

Let $d \in \{2,3\}$ and let \(W \colon \R^{d \times d} \to \R_{\infty}\) be a function which maps $(d \times d)$-matrices to real scalars or infinity. 
The notion of polyconvexity relies on the minors of the matrices \(F \in  \R^{d \times d} \). Given the determinant $\det(F)$ and the adjugate $\adj(F)$ of $F$, let 
\begin{equation} \label{eq:M(F)}
    \mathcal{M}(F)=
    \begin{cases}
    (F, \det(F)) & \text{if } \  d=2,\\
    (F, \adj(F), \det(F))  & \text{if } \ d=3      
    \end{cases}
\end{equation}
denote the minors of $F$, i.e.~$\M(F)$ is a vector of dimension $K_d = 5$ if $d =2$ and $K_d = 19$ if $d=3$.
A function \(V \colon \R^{d \times d} \to \R_{\infty}\) is said to be polyconvex if there exists a convex function 
\(G\colon \R^{K_d} \to \R_{\infty}\) such that for all \(F \in \R^{d \times d}\), 
\begin{equation*}
    V(F) = G(\mathcal{M}(F)).
\end{equation*}
The polyconvex envelope \(W^{\pc}(F) \colon \R^{d \times d} \to \R_{\infty}\) of \(W\), defined by the pointwise supremum
\begin{equation} \label{eq:defWpc}
    W^{\pc}(F) = \sup \{V(F) \ \lvert \ V \colon \R^{d \times d} \to \R_{\infty} \ \text{polyconvex}, V \leq W \},
\end{equation}
is the largest polyconvex function below $W$. 

We will restrict ourselves to the class of isotropic densities and call $W$ isotropic if and only if 
\begin{equation*}
	W(F) = W(R_1 F R_2)
\end{equation*}
for all $F\in \R^{d \times d}$ and all $R_1, R_2 \in \SOd$ and we say that $W$ is $\SOd \times \SOd$-invariant in this case, where $\SOd$ denotes the special orthogonal group of $(d \times d)$-matrices. 
Note that we define isotropy in the sense as introduced in \cite{Bal76}, the left \(\SOd\)-invariance includes objectivity, while the right invariance reflects full material symmetry.
Following \cite{WiePet2026}, isotropic functions can be characterised by the signed singular values of their arguments. Let $0\leq \sigma_1,\ldots,\sigma_d\in\R$ denote the singular values of a matrix $F\in\R^{d\times d}$.
Then \(\nu_1,\ldots, \nu_d \in \R\) are called \emph{signed singular values} of \(F\), if they have the same absolute values as the singular values of $F$, i.e.~up to permutation it holds \(|\nu_i| = \sigma_i\), and their signs satisfy \(\sign(\nu_1 \cdot \ldots \cdot \nu_d) = \sign(\det(F))\).
In this definition, the signed singular values are only unique up to permutations in 
\begin{equation*}
	\Pid = \left\{P \diag(\varepsilon) \in \mathcal{O}(d) \mid P \in \operatorname{Perm}(d), \varepsilon \in \{-1,1\}^d, \varepsilon_1\cdot\ldots\cdot\varepsilon_d = 1 \right\},
\end{equation*}
where $\diag$ refers to the diagonal matrix with entries given by the vector of its argument and $\operatorname{Perm}(d) \subset \{0,1\}^{d \times d}$ denotes the set of  permutation matrices. 
In what follows, we denote by \(\nu \colon \R^{d \times d} \to \R^{d}\) the signed singular value mapping assumed to be well-defined, e.g.~by fixing the order of the entries according to the magnitude of the absolute values and assuming positivity of the entries up to only one entry, i.e.~assigning the sign of the determinant to a single fixed entry.

It is possible to identify the set of isotropic functions $W\colon \R^{d \times d} \to \R_{\infty}$ with the set of $\Pid$-invariant functions $\Phi\colon\R^d \to \R_{\infty}$, i.e.~${\Phi(\hat{\nu}) = \Phi(S \hat{\nu})}$ for all $\hat{\nu}\in \R^{d}$ and all $S \in \Pid$. 
The identifications are given by
\begin{equation} \label{eq:W=PhiPhi=W}
	W(F) = \Phi(\nu(F)) \qquad \text{and} \qquad \Phi(\hat{\nu}) = W(\diag(\hat{\nu}))
\end{equation}
for all $F\in \R^{d\times d}$ and for all vectors $\hat{\nu} \in \R^{d}$. 

As shown in \cite{WiePet2026}, the polyconvexity of isotropic functions can be characterised directly using a lower dimensional mapping \(\Phi\). For this purpose, in analogy to the minors \(\M\) in \eqref{eq:M(F)}, we define \( k_d \coloneqq 2^d - 1\), and introduce the mapping \(\m\colon \R^{d} \to \R^{k_d}\) by 
\begin{equation*}
	\m(\hat{\nu}) = 
	\begin{cases}
		({\nu}_1, {\nu}_2, {\nu}_1 \, {\nu}_2) & \text{if } \ d = 2, \\
		({\nu}_1, {\nu}_2, {\nu}_3, {\nu}_2 \, {\nu}_3, {\nu}_3 \, {\nu}_1, {\nu}_1 \, {\nu}_2, {\nu}_1 \, {\nu}_2 \, {\nu}_3) & \text{if } \ d = 3.
	\end{cases}
\end{equation*}
We call $\m(\hat{\nu})$ the vector of minors of $\hat{\nu} = (\nu_1,\ldots, \nu_d)\in\R^d$; in the literature, this vector is also called the elementary polynomials. 
A \(\Pid\)-invariant function \(\Psi\) is called (signed singular value) polyconvex if there exists a convex function \(g\colon \R^{k_d} \to \R_{\infty}\) such that \(\Psi = g\circ\m\), see \cite[Corollary~1.5]{WiePet2026}.
The (signed singular value) polyconvex envelope of \(\Phi\) is then defined by
\begin{equation}\label{eq:defPhipc}
    \Phi^{\pc}(\hat{\nu}) = \sup\{\Psi(\hat{\nu}) \mid \Psi \colon \R^{d} \to \R_{\infty} \text{ polyconvex }, \Psi \leq \Phi\}.
\end{equation}
According to \cite[Remark~2.2]{NPPW2024}, the polyconvex envelope $W^{\pc}$ of the original density can be identified with $\Phi^{\pc}$, similar to \eqref{eq:W=PhiPhi=W}, in the following sense:
\begin{equation} \label{eq:Wpc=PhipcPhipc=Wpc}
	W^{\pc}(F) = \Phi^{\pc}(\nu(F)) \qquad \text{and} \qquad \Phi^{\pc}(\hat{\nu}) = W^{\pc}(\diag(\hat{\nu})).
\end{equation}
So we can limit ourselves to the approximation of \(\Phi^{\pc}\). Thanks to  \cite[Corollary~2.5]{NPPW2024}, it can be obtained by
\begin{equation} \label{eq:Phipc=hc}
    \Phi^{\pc}(\hat{\nu}) = h^{\cc} (\m (\hat{\nu})),
\end{equation}
where \(h^{\cc}\) is the convex envelope of the function
\begin{equation*}
    h\colon \R^{k_d} \to \R_\infty,
    \qquad 
    x \mapsto
    \begin{cases}
        \Phi(\hat{\nu}) &\text{if } \ x = m(\hat{\nu}), \\
		\infty &\text{else}.
    \end{cases}
\end{equation*}
The relation \eqref{eq:Phipc=hc} transforms the polyconvexification problem of \(\Phi\) into a convexification problem of the function \(h\) in the lifted signed singular value space. Compared to the original polyconvexification problem defined on \(d \times d\) matrices, this formulation reduces the dimension of the manifold to be convexified from \(d \times d\) to \(d\). In particular, as shown in \cite{NPPW2024}, the dimension of the ambient space is reduced from \(19\) and \(5\) to \(7\) and \(3\) for \(d=3\) and \(d=2\), respectively, allowing an efficient numerical treatment.

\section{Properties-Preserving Neural Networks} \label{sec:theoryNN}
Given the notion of polyconvexity for isotropic functions introduced in the previous section, we want to approximate the parameter-dependent polyconvex envelope \(\Phi^{\pc}\colon \R^{d}\times \R^{p} \to \R\) using a neural network denoted by \(\Phi^{\pc}_{\pred}\colon \R^{d}\times \R^{p} \to \R\). This is achieved by approximating the convex envelope \(h^{\cc}\) of the function \(h\colon \R^{k_d}\times \R^{p} \to \R_{\infty}\) defined in \eqref{eq:Phipc=hc}. For clarity, we restrict ourselves to the case where both \(\Phi^{\pc}\) and \(\Phi^{\pc}_{\pred}\) are finite and do not attain the value infinity. In some expressions, we simplify the notation by writing \(\hat{m} \coloneqq \m(\hat{\nu})\) for the minors of the signed singular values vector \(\hat{\nu}\).

\subsection{Neural Network Concepts}
An artificial neural network is a mathematical model defined as a composition of functions, each of which is typically formed by further compositions of functions \cite{BookNN}. This model can be conveniently represented as a network structure. In what follows, we present the basic notions of layers, weights, inputs and outputs of a general neural network.

Consider a neural network consisting of \(k\) layers. 
We denote the index of layers by the superscript \(i \in \{0,\ldots,k\}\), where \(i = 0\) denotes the input layer, \(i = 1\) the first hidden layer and \(i = k\) the output layer. 
The total number of hidden layers is thus \(k-1\). 
The number of neurons in the \(i\)-th layer is denoted by \(d_{i}\); in particular, \(d_{0}\) is the number of inputs and \(d_{k}\) is the number of outputs.
Let \(W_{i} \in \R^{d_{i+1} \times d_{i}}\) denote the weight matrix associated with the connection from the \(i\)-th to the \((i+1)\)-th layer, \(b_{i} \in \R^{d_{i + 1}}\) the corresponding bias vector and \(g_i\) the activation function, acting component-wisely. 
For \(i = 0, \dots, k-1\), the network recurrence on the layer outputs \(z_i \in \R^{d_{i}}\) can be written in matrix–vector form as
\begin{equation*}
	z_{i + 1} = g_{i}\!\left(W_{i} \, z_{i} + b_{i} \right),
\end{equation*}
where \(z_{0}\) denotes the input and \(z_{k}\) the overall output of the network.

Building on the characterisation of the (signed singular value) polyconvex envelope \(\Phi^{\pc}\) in \cref{sec:theorypoly}, we aim to construct a neural network approximation \(\Phi^{\pc}_{\pred}\) that reliably preserves its physical properties. The architecture must ensure polyconvexity, i.e.~convexity with respect to the minors of the signed singular values, and satisfy the envelope inequality $\Phi^{\pc} \leq \Phi$ (pointwisely) as in \eqref{eq:defPhipc}. In addition, the network must respect the \(\Pi_{d}\)- invariance of \(\Phi\) representing symmetries. The following sections outline strategies for enforcing these three key properties into the design of neural networks.

\subsection{Enforcing Polyconvexity}
To enforce the convexity in the minors \(\m(\hat{\nu})\) of the input vector \(\hat{\nu}\), a particular class of neural networks is employed: Input Convex Neural Networks (ICNN) introduced in \cite{ICCN}. 
In what follows, two variants of ICNN are considered: the Fully Input Convex Neural Networks (FICNN) and the Partially Input Convex Neural Networks (PICNN). 

\subsubsection{Fully Input Convex Neural Networks (FICNN)}
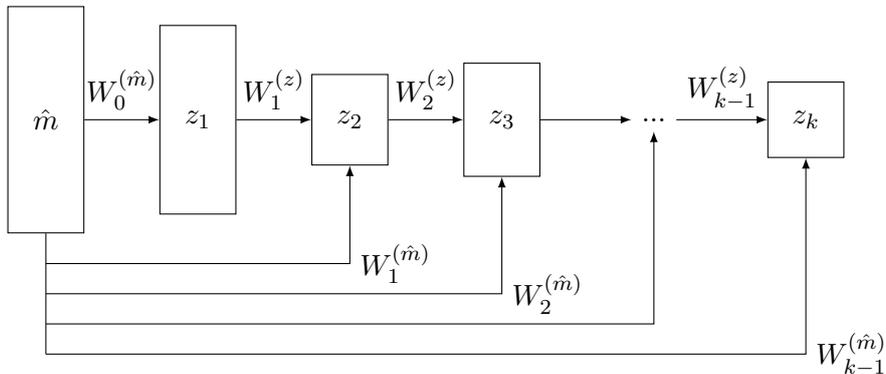
\begin{figure}[htbp]
    \begin{center}
    \begin{tikzpicture}[>=latex, scale =1]
    \node[draw, minimum width=1cm, minimum height=3cm, align=center] (y) at (0, 0) {$\hat{m}$};
    \node[draw, minimum width=1cm, minimum height=2.5cm, align=center] (z1) at (2, 0) {$z_1$};
    \node[draw, minimum width=1cm, minimum height=1.2cm, align=center] (z2) at (4, 0) {$z_2$};
     \node[draw, minimum width=1cm, minimum height=1.5cm, align=center] (z3) at (6, 0) {$z_3$};
     \node[align=center] (zq) at (8, 0) {$...$};
    \node[draw, minimum width=1cm, minimum height=1cm, align=center] (zk) at (10, 0) {$z_k$};

    \draw[->] (y) -- node[above] {$W_0^{(\hat{m})}$} (z1);
    \draw[->] (z1) -- node[above] {$W_1^{(z)}$} (z2);
    \draw[->] (z2) -- node[above] {$W_2^{(z)}$} (z3);
    \draw[->] (z3) -- node[above] {} (zq);
    \draw[->] (zq) -- node[above] {$W_{k-1}^{(z)}$} (zk);
   
    \draw[->] (y.south) ++(0, -0.4) -| node[below, right] {$W_1^{(\hat{m})}$} (z2.south);
    \draw[->] (y.south) ++(0, -0.8) -| node[below, right] {$W_2^{(\hat{m})}$} (z3.south);
    \draw[->] (y.south) ++(0, -1.2) -| node[below, right] {}(zq.south);
    \draw[->] (y.south) ++(0, -1.6) -| node[below, right] {$W_{k-1}^{(\hat{m})}$}(zk.south);
    \draw[-] (y.south) ++(0, 0) -- ++(0, -1.6);
 
\end{tikzpicture}
    \end{center}
    \caption{A Fully Input Convex Neural Network (FICNN).}
    \label{fig:FICNN}
\end{figure}

The model illustrated in \cref{fig:FICNN} defines a neural network over the input $\hat{m}$ using the recurrence for $i=0, \dotsc, k-1$, 
\begin{equation*}
    z_{i+1} = g_i(W_i^{(z)} z_i + W_i^{(\hat{m})} \hat{m} + b_i), 
\end{equation*}
where $z_i$ denotes the layer activations (with $z_0 = 0$ , $W_0^{(z)} \equiv 0$) and $g_i$ are nonlinear activation functions.
The overall network evaluation reads 
\begin{equation*}
    \Phi^{\pc}_{\pred}(\hat{\nu}) \coloneqq h^{\rm c}_\pred(\hat{m};\theta) = z_k,
\end{equation*}
where $h^{\rm c}_\pred$ is the neural network approximation of the function $h^{\rm c}$ defined in \eqref{eq:Phipc=hc} and depends on $\theta = \{ W_{0:k-1}^{(\hat{m})}, W_{1:k-1}^{(z)}, b_{0:k-1} \}$, the parameter vector collecting the weights and biases. 

The convexity of $h^{\cc}_\pred$ is addressed first. In order for $\Phi^{\pc}_{\pred}$ to be polyconvex as defined in \cref{sec:theorypoly}, the $\Pid$-invariance is also required, which will be addressed in \cref{sec:pid}.

\begin{prop}
    The function 
    $h^{\cc}_\pred(\,\cdot\,;\theta)$ is convex in $\hat{m}$, i.e.~the minors $\m(\hat{\nu})$ of the signed singular values vector \(\hat{\nu}\), provided that all $W_{1:k-1}^{(z)}$ are non-negative and all activation functions $g_i$ are convex and non-decreasing.
\end{prop}

The proof follows from the fact that non-negative linear combinations of convex functions are also convex, and that the composition of a convex and convex non-decreasing function is also convex.  The constraint that the $g_i$ are convex and non-decreasing is not particularly restrictive, since common nonlinear activation functions such as the ReLU or SoftPlus activation functions satisfy this constraint. A review of convex activation functions is presented in \cite{Wilhelm}.

\subsubsection{Partially Input Convex Neural Networks (PICNN)}
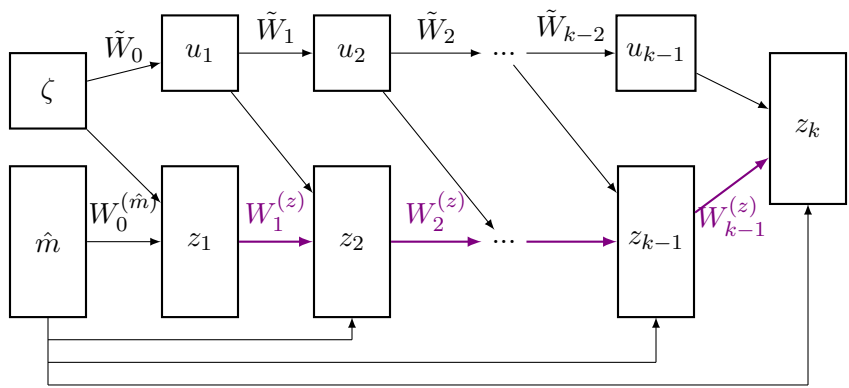
\begin{figure}[htbp]
    \begin{center}
    \begin{tikzpicture}[>=latex, scale=1]

\node[draw, thick, minimum width=1cm, minimum height=1cm, align=center] (I2)    at (0,2) {$\zeta$};

\node[draw, thick, minimum width=1cm, minimum height=2cm, align=center] (I1) at (0,0) {$\hat{m}$};

\node[draw, thick, minimum width=1cm, minimum height=1cm, align=center] (G1) at (2,2.5) {$u_1$};

\node[draw, thick, minimum width=1cm, minimum height=1cm, align=center] (G2) at (4,2.5)  {$u_2$};

\node[align=center] (G3) at (6,2.5) {$...$};

\node[draw, thick, minimum width=1cm, minimum height=1cm, align=center] (G4) at (8,2.5)  {$u_{k-1}$};

\node[draw, thick, minimum width=1cm, minimum height=2cm, align=center] (J1) at (2,0) {$z_1$};

\node[draw, thick, minimum width=1cm, minimum height=2cm, align=center] (J2) at (4,0) {$z_2$};

\node[align=center] (J3) at (6,0) {$...$};

\node[draw, thick, minimum width=1cm, minimum height=2cm, align=center] (J4) at (8,0) {$z_{k-1}$};

\node[draw, thick, minimum width=1cm, minimum height=2cm, align=center] (J5) at (10,1.5) {$z_k$};

\draw[->] (I2) -- node[above] {$\tilde{W}_0$} (G1);
\draw[->] (G1) -- node[above] {$\tilde{W}_1$} (G2);
\draw[->] (G2) -- node[above] {$\tilde{W}_2$} (G3);
\draw[->] (G3) -- node[above] {$\tilde{W}_{k-2}$} (G4);
\draw[->] (G4) -- (J5);

\draw[->] (I1) -- node[above] {$W_0^{(\hat{m})}$} (J1);
\draw[thick, violet, ->] (J1) --  node[above] {$W_1^{(z)}$} (J2);
\draw[thick, violet,->] (J2) --  node[above] {$W_2^{(z)}$} (J3);
\draw[thick, violet,->] (J3) -- (J4);
\draw[thick, violet,->] (J4) --  node[below] {$W_{k-1}^{(z)}$} (J5);

\draw[->] (I2) --   (J1);
\draw[->] (G1) -- (J2);
\draw[->] (G2) -- (J3);
\draw[->] (G3) -- (J4);

\draw (0,-1) -- (0,-1.9);
\draw (0,-1.3) --   (4,-1.3);
\draw (0,-1.6) --  (8,-1.6);
\draw (0,-1.9) --   (10,-1.9);

\draw[->] (4,-1.3) -- (4,-1);
\draw[->] (8,-1.6) -- (8,-1);
\draw[->] (10,-1.9) -- (10,0.5);

\end{tikzpicture}
    \end{center}
    \caption{A Partially Input Convex Neural Network (PICNN), coloured weights need to be non-negative to ensure convexity of the network in \(\hat{m}\).}
    \label{fig:PICNN}
\end{figure}

The previous model, i.e.~FICNN, provides convexity over the entire input of the neural network, which may in fact be a restriction on the allowable class of models. 
Furthermore, this full joint convexity is unnecessary in the setting where the convexity is required only over some inputs.  
This is why we present the Partially Input Convex Neural Networks (PICNN), as introduced in \cite{ICCN}, which are convex over only some inputs of the network. The $k$-layer PICNN architecture (\cref{fig:PICNN}) is defined by the recurrence for $i=0, \dotsc, k-1$,
\begin{equation*}
\begin{aligned}
    &u_{i+1} = \tilde{g}_i \left (\tilde{W}_i u_i + \tilde{b}_i \right), \\
    &z_{i+1} = g_i\left( W_i^{(z)}\Big(z_i \circ  \max([W_i^{(zu)}u_i+b_i^{(z)}],0)\Big) +  W_i^{(\hat{m})}
    \Big(\hat{m} \circ [W_i^{(\hat{m}u)}u_i+b_i^{(\hat{m})}]\Big) + W_i^{(u)}u_i + b_i \right), 
\end{aligned}
\end{equation*}
with $u_0 =\zeta$, $z_0 =0$ and $W_0^{(z)}=0$ and
\begin{equation*}
    \Phi^\pc_\pred(\hat{\nu}; \zeta) \coloneqq h_\pred^{\rm c}(\hat{m};\zeta;\theta) = z_k
\end{equation*}
denotes the overall evaluation of the network. Moreover, $h_\pred^{\cc}$ is the neural network approximation of the function $h^{\cc}$ defined in \eqref{eq:Phipc=hc} and depends on $\theta$, the parameter vector collecting the weights and biases. The $u_i \in \R^{n_i}$ and $z_i \in \R^{m_i}$ denote the hidden units for the $\zeta$-path and $\hat{m}$-path, and $\circ$ denotes the Hadamard product, i.e.~the element-wise product between two vectors. 
The functions $g_i$ and $\tilde{g}_i$ are activation functions.

\begin{prop}
    The function
    $h_{\pred}^{\cc}(\,\cdot\,;\zeta;\theta)$ is convex in $\hat{m}$, i.e.~the minors $m(\hat{\nu})$ of the signed singular values vector \(\hat{\nu}\), provided that the weights $W_{1:k-1}^{(z)}$ are non-negative and all activation functions $g_i$ associated with the \(\hat{\m}\) path are convex and non-decreasing. 
\end{prop}

Note that the schemes shown in \cref{fig:FICNN,fig:PICNN} already reflect the convexity with respect to the minors \(\hat{\m}\), i.e.~the formulation to ensure polyconvexity. However, the change from the original input \((\hat{\nu}, \zeta)\) of the function \(\Phi^\pc\) to the input \((\m(\hat{\nu}), \zeta)\) of the function \(h^{\cc}\), i.e.~the transitioning from the signed singular values to the minors lifted signed singular values, can be interpreted as a hard-coded feature extracting layer in the overall network architecture. \\

\begin{rem}
  \cite{chen2019,huang2021} prove universal approximation theory of PICNNs when ReLU, as well as SoftPlus, activation functions are used. Besides, \cite{GeuKurWieMos2025} provide a proof of a universal approximation for frame-indifferent isotropic polyconvex functions for ICNNs combined with signed singular values. These results are directly applicable to the architecture considered in this work, which also justifies our architectural choice.
\end{rem}

In the present work, FICNN is employed when the function \(\Phi^{\pc}\) to be predicted does not depend on parameters but only on the signed singular values vector $\hat{\nu}$, and PICNN is employed when the function to be predicted depends on the signed singular values vector $\hat{\nu}$ along with additional parameters \(\zeta\). 
Indeed, the use of PICNN allows ensuring convexity with respect to the minors $m(\hat{\nu})$ of the signed singular values vector $\hat{\nu}$ while relaxing the constraints on the $\zeta$-path, allowing a better representation.

\begin{rem}
Besides Input Convex Neural Networks, several alternative strategies exist to enforce convexity or polyconvexity. One can add a Hessian‐based penalty to the loss to impose positive‐definiteness and hence convexity \cite{Vahidullah}. However, it is important to note that polyconvexity is a global property and therefore cannot be directly enforced through a loss term. 
Enforcing it would require relaxing the condition to a weaker notion, such as local polyconvexity, see \cite{KalGebBruLinSunKaes2024}.
\cite{Vahidullah2} implements a neural‐ODE formulation that embeds polyconvexity directly into the model structure so that it holds automatically. Another approach defines an integral operator that maps any base function into a strongly convex surrogate a priori \cite{chen2022}. Finally, the use of polyconvex activation functions combined with tailored network architectures guarantees polyconvex strain‐energy functions by construction \cite{Linka}.
\end{rem}

\subsection{Enforcing the Envelope Inequality} \label{sec:PenUpper}
By definition \eqref{eq:defPhipc}, the polyconvex envelope $\Phi^\pc$ of the function $\Phi$ is the largest polyconvex function below $\Phi$. Consequently, the goal is to ensure that the polyconvex envelopes predicted by the neural network are below the function \(\Phi\), i.e.~\(\Phi^{\pc}_{\pred} \leq \Phi\) which we refer to as envelope inequality. 
In the present work, the envelope inequality is imposed in a weak sense by penalisation, avoiding the addition of constraints directly in the architecture. 
To achieve this, a custom loss function is designed, which penalises predictions of the neural network located above the function $\Phi$ during the training process.
By minimising this loss function, the network ultimately achieves predictions below \(\Phi\). The loss function is thus written as 
\begin{equation} \label{eq:loss1}
    \Loss = \Loss_{\mse} + \lambda_{\ineq} \, \Loss_{\ineq}
\end{equation}
with 
\begin{equation}
    \Loss_{\mse} = \frac{1}{N} \sum_{i=1}^N (\Phi^{\pc}_{i} - \Phi^{\pc}_{\pred,i})^2 ,
\label{eq:Lmse}
\end{equation} 
the classical mean square error, and 
\begin{equation} \label{eq:lossineq}
    \Loss_{\ineq} = \frac{1}{N} \sum_{i=1}^N \max \{\Phi^{\pc}_{\pred,i} - \Phi_i, 0 \}^2 ,
\end{equation}
where $N$ is the size of the learning data, $\Phi^{\pc}_{\pred}$ is the neural network prediction of the polyconvex envelope, and $\Phi^{\pc}$ is the target value to be approximated.
The subscript \(i\) denotes the evaluation of these functions at the \(i\)th learning input data point given as a tuple of the form \((\m(\hat{\nu}), \zeta)\). 
The parameter \(\lambda_{\ineq} \geq 0\) is the so-called penalty parameter, which determines how the two terms are weighted against each other and determines how strongly the network should be penalised if the prediction lies above \(\Phi\).  

\begin{rem} \label{rem:alternativesEnvelopeIneq}
    Alternatively to the penalisation via \eqref{eq:lossineq}, the envelope inequality can be enforced by embedding constraints directly into the network architecture.
    A multiplication layer whose output lies in the interval \([0, 1]\) can be employed to ensure the inequality via multiplication with the function value of \(\Phi\).
    Typically, the required activation function, ensuring outputs in \([0, 1]\), is non convex and hence might conflict with the polyconvexity.
    Another viable way is to add a minimum operation layer, i.e.~the final output is derived by computation of the pointwise minimum between the prediction and the function value, i.e.~\(\min\{\Phi^{\pc}_{\pred}, \ \Phi\}\). 
    However, numerical experiments have shown that this constraint is too restrictive and the model stops learning if the prediction is already below the function \(\Phi\).
\end{rem}

\subsection{Enforcing the \(\Pi_{d}\)-Symmetry}\label{sec:pid}
The function \(\Phi\) from \eqref{eq:W=PhiPhi=W} is subject to the invariance of \(\Pi_{d}\), and the same holds for \(\Phi^{\pc}\) from \eqref{eq:Wpc=PhipcPhipc=Wpc}, that is it holds
\begin{equation*}
    \Phi^{\pc}(\pi(\hat{\nu})) = \Phi^{\pc}(\hat{\nu})
\end{equation*}
for all \(\pi \in \Pi_{d}\).
Consequently, this symmetry should be preserved by a neural network approximation to ensure the \(\Pid\)-invariance as required for the notion of (signed singular value) polyconvexity, as defined in \cref{sec:theorypoly}.

As for the envelope inequality in \cref{sec:PenUpper}, the symmetry of the prediction can be imposed in a weak sense by penalisation.  To achieve this, a custom loss function which penalises more severely if the prediction of the neural network is not symmetric during the training process is designed. By minimising this loss function, the network should ultimately achieve symmetric predictions.  Building upon the previously defined loss function \eqref{eq:loss1} (keeping the same notations), already facilitating the envelope inequality, we add a term for the penalisation of the non-symmetry, leading to the loss function defined by
\begin{equation} \label{eq:loss}
   \Loss = \Loss_{\mse} + \lambda_{\ineq} \, \Loss_{\ineq} +  \lambda_{\sym} \, \Loss_{\sym}
\end{equation}
with
\begin{equation} \label{eq:losssym} 
   \Loss_{\sym} = \frac{1}{\lvert \Pi_d \rvert} \sum_{\pi \in \Pid} \frac{1}{N} \sum_{i=1}^{N} \left(\Phi^{\pc}_{\pred,i}(\hat{\nu}) -\Phi^{\pc}_{\pred,i}(\pi(\hat{\nu}))\right)^2,
\end{equation}
where the normalisation constant is motivated by the fact that the term corresponding to $\pi = \id$ does not contribute to the overall loss. 
The parameter $\lambda_{\sym} \geq 0$ is the penalty parameter, which determines how strongly the network should be penalised if the network is not symmetric with respect to the inputs. 
In addition, symmetry is further fostered by data augmentation. 
We incorporate symmetric data in the dataset, i.e.~for a given point $\hat{\nu}$, the points $\pi(\hat{\nu})$ for \(\pi \in \Pid\) are also included in the dataset.

\begin{rem}[Other approaches to ensure symmetry]
Another approach to ensure symmetry it to hard-code it directly into the neural network architecture, for example by employing the weight-sharing strategy as described in \cite[Appendix~B]{Vijayakumaran}. In the case $d=2$, this method  consists in using identical weights for the direct connections to the inputs $\nu_1$ and $\nu_2$, i.e.~to the signed singular values. This ensures that the network treats symmetric inputs in an equivalent manner, thereby preserving the desired symmetry. However, such an approach would ensure symmetry only for permutations, with possibly some loss of  approximation capacity, but not for symmetry which also include an even number of sign changes, such as the \(\Pid\)-symmetry. In addition, the symmetry can also be enforced a posteriori, see for example \cite{GeuKurWieMos2025}, in which the symmetry is realized by averaging over all input permutations in the symmetry group 
\begin{equation*}
    \bar{\Phi}^{\pc}_\pred(\hat{\nu}) =\frac{1}{\lvert \Pi_d \rvert} \sum_{\pi \in \Pid} \Phi^{\pc}_\pred(\pi(\hat{\nu})).
\end{equation*}
This method has the advantage of ensuring output symmetry in \(\bar{\Phi}^{\pc}_{\pred}\) without modifying the neural network architecture.
However, this requires multiple evaluations of \(\Phi^{\pc}_{\pred}\) and the approximation \(\Phi^{\pc}_{\pred}\) is not necessarily symmetric.
We do not pursue these symmetry-incorporation variants further, as the numerical experiments in \cref{sec:KSD} indicate that weak incorporation into the loss function already yields good accuracy.
\end{rem}

\subsection{Implementation and Hyperparameters}
The neural network architectures are implemented using PyTorch, and all the trainings are performed on an Intel\textsuperscript{\textregistered} Core\textsuperscript{TM} i9-11900K machine, using only one core at a base frequency of 3.50 GHz. 

The input data of the neural networks consist of the minors \(\m(\hat{\nu}) \in \R^{k_d}\) of the signed singular values \(\hat{\nu}\) along with additional parameters \(\zeta \in \R^{p}\) depending on the definition of the function \(\Phi\), i.e.~the input data consists of vectors \((\m(\hat{\nu}), \zeta) \in \R^{k_d + p}\).
The target values correspond to the evaluation of the polyconvex envelope \(\Phi^{\pc}\) at the points $\hat{\nu}$ and the parameters $\zeta$, expressed as $\Phi^{\pc}(\hat{\nu}; \zeta)$. 
In summary, the learning data (gathering both training and validation data) consists of tuples of the form \((m(\hat{\nu}), \zeta, \Phi^{\pc}(\hat{\nu}; \zeta), \Phi(\hat{\nu}; \zeta))\) containing the input features, the target value $\Phi^{\pc}(\hat{\nu}; \zeta)$, and the value of the function $\Phi(\hat{\nu}; \zeta)$, which is necessary for the evaluation of the loss function \eqref{eq:loss}. 

In all the numerical examples, we design the neural network architectures to be as small as possible while maintaining high predictive accuracy. 
Indeed, preliminary numerical experiments have shown that increasing the number of layers or neurons of the neural networks presented below does not necessarily improve prediction performance; on the contrary, it may even degrade them. 
Also, we aim to use the smallest dataset size that ensures accurate predictions, while making training computationally efficient and feasible on a standard laptop, which has been achieved thanks to several investigations. 
We employ the ReLU activation function in all hidden layers, as it satisfies the convexity condition required by ICNN. 
The output layer consists of a single neuron with a linear activation function. 
Numerical experiments indicate that selecting ReLU as the convex activation function enhances learning efficiency and improves predictive accuracy compared to the SoftPlus function, another commonly employed activation function in ICNNs.

The weights $W^{(z)}$ which belong to the convex part of the network are initialised from a normal distribution $\mathcal{N}(0.1,0.1)$ with $0.1$ mean and \(0.1\) standard deviation, and then projected onto $\R_+$.
Also, after each training step, the weights $W^{(z)}$ are projected onto the $\R_+$ halfspace to ensure their positiveness.  This projection is essential to ensure the positivity of the weights. The naive application of ReLU for weight clipping can result in zero entries in the weight matrices and would waste approximation potential. 
For this projection, a shifted version of ReLu, i.e.~\(x \to \max(x, 0) + \varepsilon \) with the small offset \(\varepsilon = 10^{-6}\) is employed.
This initialisation and projection strategy has shown better performance compared to a uniform distribution and projections using exponential or SoftPlus functions.
The other weights are initialised from a uniform distribution \(\mathcal{U}(- 1/\sqrt{n}, 1/\sqrt{n})\) with $n$ the input size and all the biases are initialised according to \(\mathcal{N}(0, 0.1)\).

We employ the Adam optimiser with a learning rate of \(\eta = 0.001\), a batch size of \(128\), with the data shuffled at each epoch to improve model generalisation, and a patience of \(5\), i.e.~the training process is stopped if no improvement of the validation loss is achieved for \(5\) successive epochs.  
Although more complex training strategies, e.g.~with the use of learning rate decay, could have been implemented, numerical experiments demonstrate that such an approach seems unnecessary in this context. The training and validation losses are computed with the loss function \eqref{eq:loss}. The hyperparameters $\lambda_{\ineq}$ and $\lambda_{\sym}$ involved in \(\Loss\) have been determined by performing several preliminary tests with different parameter values for $\lambda_{\ineq}$ and $\lambda_{\sym}$.  
In each numerical example presented below, the given values of the $(\lambda_{\ineq}, \lambda_{\sym})$ achieved favourable training performance and predictions. 

To obtain more reliable numerical results, for each example, we conduct ten or twenty independent runs of training, using the same training and validation data.  
The variations arise solely from the initialisation of model weights and the stochastic nature of data loading during training, since a shuffle is used. The numerical results presented are obtained by averaging the outputs of the ten network realisations.
To assess the accuracy of the neural network predictions \(\Phi^{\pc}_{\pred}\) in comparison to the ground truth \(\Phi^{\pc}\), obtained either analytically or by numerical algorithms, we compute the following error metrics, namely the mean error
\begin{equation*}
	{\rm Mean \ err} = 
	\frac{1}{N}\sum_{i=1}^N \lvert\Phi^{\pc}_{\pred,i}-\Phi^{\pc}_{i}\rvert ,
\end{equation*}
the relative quadratic error
\begin{equation*}
	{\rm Rel\ quad \ err} = 
	\sqrt{\frac{\sum_{i=1}^N \lvert \Phi^{\pc}_{\pred,i}-\Phi^{\pc}_{i}\rvert^2}{\sum_{i=1}^N  \lvert\Phi^{\pc}_{i}\rvert^2}} ,
\end{equation*}
and the relative maximum error
\begin{equation*}
	{\rm Rel\ max\ err} = 
	\frac{\max_{i=1,\ldots, N} \lvert \Phi^{\pc}_{\pred,i}-\Phi^{\pc}_{i} \rvert}{\max_{i=1,\ldots, N} \lvert\Phi^{\pc}_{i}\rvert},
\end{equation*}
where \(N\) denotes the number of evaluation points.

\subsection{Data Generation} 
\label{sec:SVPC}
For the physically relevant benchmark examples considered in the next sections, no closed-form analytical representations of the polyconvex envelopes are available. 
This is why, for the computation of the target values the polyconvex envelopes are approximated numerically utilising the algorithm for isotropic functions based on a linear programming approach presented in \cite{NPPW2024}. 
For this approximation procedure, the pointwise characterisation of the polyconvex envelope \(\Phi^{\pc}\) at \(\hat{\nu} \in \R^{d}\) by the optimisation problem
\begin{equation} \label{eq:poly-opt-prob-iso}
	\Phi^{\pc}(\hat\nu) = \inf \left\{\sum_{i = 1}^{k_d + 1} \xi_i \, \Phi(\nu_i)\,\biggl\vert\, \xi_{i} \in [0, 1],\, \nu_i \in \R^{d},\, \sum_{i = 1}^{k_d + 1} \xi_{i} = 1,\, \sum_{i = 1}^{k_d + 1} \xi_i \, \m(\nu_i) = \m(\hat{\nu})\right\}
\end{equation}
is employed. 
For an algorithmic approximation of this problem, consider a discretisation of the signed singular value space by a point cloud, denoted by \(\Sigma_{\delta} = \{\nu_1, \ldots, \nu_{N_\delta}\} \subset \R^{d}\), a possible choice is a structured lattice on \([-r, r]^{d}\) with lattice size \(\delta\) and discretisation radius \(r\). 
The lifting of this point cloud \(\m(\Sigma_\delta)\) can be employed to turn the nonlinear optimisation problem \eqref{eq:poly-opt-prob-iso} into the following linear program
\begin{equation}\label{eq:linearprogram}
	\Phi^{\pc}_{\delta} (\hat\nu) = \min \left\{\sum_{i = 1}^{N_\delta} \xi_{i} \, \Phi(\nu_i) \;\biggl\vert\; \xi_{i} \geq 0,\, \sum_{i = 1}^{N_\delta} \xi_{i} = 1,\, \sum_{i = 1}^{N_\delta} \xi_{i}\,\m(\nu_{i}) = \m(\hat\nu) \right\},
\end{equation}
which can be efficiently solved numerically using standard algorithms for linear programming.
The overall procedure is referred to as signed singular value polyconvexification by linear programming (SVPC~LP). 
In exact arithmetic, there exists a minimiser \(\xi \in \R^{N_\delta}\) of \eqref{eq:linearprogram} with at most \(k_d+1\) non-zero entries reflecting the convex coefficients of the supporting points of the polyconvex envelope at \(\hat{\nu}\) in the set \(\Sigma_{\delta}\) which are sometimes also referred to as volume fractions.
For a detailed description of the algorithm see \cite{NPPW2024}, a MATLAB and Python implementation can be found under \url{https://github.com/TmNmr/SVPC}.

In this paper, we use a python implementation of this algorithm which uses the \texttt{linprog} function from \texttt{scipy.optimize} for solving the linear program. 
As discretisation of the signed singular value space, the shifted lattice generated by \(\Sigma_{\delta} = (\delta\, \Z^{d} + \tfrac{\delta}{2} \, \mathbbm{1}_{d}) \cap [-r, r]^d\) is employed. 
The lattice parameters \(r\) and \(\delta\) were chosen such that an absolute error of order \(10^{-5}\) can be predicted, see e.g.~\cite[Figures~4.2,~4.4 and~4.5]{NPPW2024}. Moreover the symmetry induced by the \(\Pid\)-invariance is exploited, reducing the computation to grid points \(\hat{\nu} \in \Sigma_{\delta}\) which are located in the signed singular value cone \(0 \leq \lvert \nu_1 \rvert \leq \nu_2 \leq \ldots \leq \nu_d \leq r\) and exploitation of determinant constraints reduces the evaluation domain further to \(0 \leq \nu_1 \leq \nu_2 \leq \ldots \leq \nu_d \leq r\).
It should be noted that for points close to the coordinate axes, a sufficiently fine resolution of the signed singular value space is required due to the determinant constraints, that is why finer lattice widths \(\delta\) are employed for these regions. 
Nevertheless, in what follows, we denote by \(\Phi^{\pc}_{\delta}\) simply the approximation of the polyconvex envelope by the SVPC~LP algorithm. In addition, since SVPC~LP delivers independent point approximations, the computation of the learning data is performed in parallel. 

\section{Mathematical Benchmark Examples} \label{sec:numericalexperimentsmath}
To show the potential of the approximation approaches of the property-preserving neural networks as introduced \cref{sec:theoryNN}, we consider well-known examples in the context of relaxation algorithms in both two and three spatial dimensions. 
We consider experiments of increasing complexity: parameter-independent cases in two and three spatial dimensions, followed by parameter-dependent problems. 
Throughout all of the numerical experiments, we denote the analytical polyconvex envelope by \(\Phi^{\pc}\) and by \(\Phi^{\pc}_{\pred}\) the neural network prediction.

\subsection{The Kohn--Strang--Dolzmann Function}
\label{sec:KSD}
This function was studied in \cite{KohStr86a, KohStr86b}, subsequently modified to achieve continuity in \cite{Dol99, DolWal00} and further studied in \cite{Bar05} and serves as a first benchmark example to illustrate the capability of the neural network approach to approximate the polyconvex envelope.
We consider the function $W \colon \R^{2\times 2} \to \R$ , defined as
\begin{equation*}
    W(F)=
    \begin{cases}
        1+\lvert F \rvert^2 & \text{if } \, \lvert F \rvert \geq \sqrt{2}-1,\\
        2\sqrt{2} \, \lvert F \rvert & \text{else},
    \end{cases}
\end{equation*}
where $\lvert F \rvert \coloneqq (\sum_{i,j=1}^{d} F_{ij}^2)^{1/2}$ denotes the Frobenius norm of the matrix $F$. 
The polyconvex envelope of $W$ is explicitly known, see e.g.~\cite{Dol99}, and reads
\begin{equation*}
    W^{\pc}(F)=
    \begin{cases}
        1+\lvert F \rvert^2 & \text{if } \, \rho(F) \geq 1, \\ 
        2 \,(\rho(F) - \lvert {\det}(F) \rvert)  & \text{else},
    \end{cases}
\end{equation*}
where $\rho(F) \coloneqq\sqrt{\lvert F \rvert^2 + 2 \lvert {\det}(F) \rvert}$. 
The functions $W$ and $W^{\pc}$ are isotropic; consequently, they can be rewritten in terms of the signed singular values and reduce to $\Phi$, $\Phi^{\pc} \colon \R^2 \to \R$ with 
\begin{equation*}
    \Phi(\hat{\nu}) =
    \begin{cases}
        1+ \nu_1^2 +\nu_2^2 & \text{if } \,  \sqrt{\nu_1^2 +\nu_2^2} \geq \sqrt{2}-1,\\
        2\sqrt{2} \, \sqrt{\nu_1^2 +\nu_2^2}   & \text{else},
    \end{cases}
\end{equation*}
and 
\begin{equation}
    \Phi^{\pc}(\hat{\nu})=
    \begin{cases}
        1+ \nu_1^2 +\nu_2^2 & \text{if } \, \rho(\hat{\nu}) \geq 1, \\
        2 \, (\lvert \nu_1 \rvert + \lvert \nu_2 \rvert - \lvert \nu_1 \nu_2 \rvert) & \text{else},
    \end{cases}
\label{eq:PhipcKSD}
\end{equation}
where $\rho(\hat{\nu}) = \lvert \nu_1 \rvert + \lvert \nu_2 \rvert$.

\subsubsection{Network Architecture and Learning Data}
In this example, we implement a FICNN which consists of two hidden layers, with 10 and 20 neurons, respectively, as presented in \cref{fig:Architecture_Loss_KSD}. 
The learning domain for $\hat{\nu}$ is defined as the box $[-\overline{\nu},\overline{\nu}]^2$ with $\overline{\nu}=1.05$. For the training, the set of signed singular values \( \hat{\nu} \) used as input is generated by discretising each axis with $751$ points, so that the point 0 is included, making \(751^2\) training data points in total. 
Instead of using a uniform distribution of the points, we use a local refinement towards the origin, utilising a quadratic transformation. 
This refinement approach increases the data density near the point of non-differentiability of both functions \(\Phi\) and \(\Phi^{\pc}\). 
For the validation dataset, we randomly sample \num{169200} points within the domain $[-\overline{\nu},\overline{\nu}]^2$, accounting for $30 \%$ of the training dataset size. 
The target values, i.e.~the polyconvex envelope, are obtained by evaluation of the analytical function $\Phi^\pc$ as in \eqref{eq:PhipcKSD}.
For the loss function, we choose the penalty parameters \(\lambda_{\ineq} = 50\) and \(\lambda_{\sym} = 10\).

\subsubsection{Numerical Results}
\begin{figure}[htbp]
    \centering
    \begin{subfigure}[b]{0.49\textwidth}
    \centering
    \begin{tikzpicture}[scale=0.75, every text node part/.style={align=center}]
	\pgfmathsetmacro{\ylabels}{-3.5}
	\colorlet{colInput}{unia-pink}
	\colorlet{colHI}{unia-pink}
	\colorlet{colHII}{unia-blue}
	\colorlet{colOut}{unia-pink}

    \foreach \i in {1,2,3} {
        \node[circle, draw, fill=colInput!20, minimum size=1cm] (I\i) at (0,-1.5*\i+1.5*2) {};
    }
    \foreach \i in {1,2,3} {
        \node[circle, draw, fill=colHI!50, minimum size=0.5cm] (H\i) at (3, -\i+3.25) {};
    }
    \node at (3, -5+4.125) {\vdots};
    \node[circle, draw, fill=colHI!50, minimum size=0.5cm] (Hi) at (3, -6+3.75) {};
    
    \foreach \i in {1,2,3,4} {
        \node[circle, draw, fill=colHII!50, minimum size=0.5cm] (J\i) at (6, -\i+3.75) {};
    }
    \node at (6, -6+4.75) {\vdots};
    \node[circle, draw, fill=colHII!50, minimum size=0.5cm] (Ji) at (6, -7+4.375) {};
 
    \node[circle, draw, fill=colOut!30, minimum size=0.5cm] (O1) at (9,0) {};
      
    \foreach \i in {1,2,3} {
        \foreach \j in {1,2,3} {
            \draw[->] (I\i) -- (H\j);
        }
       	\draw[->] (I\i) -- (Hi);
    }

    \foreach \i in  {1,2,3} {
        \foreach \j in {1,2,3,4} {
            \draw[->] (H\i) -- (J\j);
        }
        \draw[->] (H\i) -- (Ji);
    }

    \foreach \j in {1,2,3,4} {
        \draw[->] (Hi) -- (J\j);
    }

    \foreach \i in  {1,2,3,4} {
        \draw[->] (J\i) -- (O1);
    }
    
    \draw[->] (Hi) -- (Ji);   
    \draw[->] (Ji) -- (O1);
    
    \node at (0,\ylabels) {Input \\ \(\hat{m}\)};

    \node at (3,\ylabels) {10 neurons \\ ReLU};

    \node at (6,\ylabels) {20 neurons \\ ReLU};

    \node at (9,\ylabels) {1 neuron \\ Linear};

    \node at (0,1.5) {$\nu_1$};
    \node at (0,0) {$\nu_2$};
    \node at (0,-1.5) {$\nu_1 \nu_2$};

    \draw[draw=none, fill=red!20, opacity=0.5] (5,-4) rectangle ++(-1,8);
    \draw[draw=none, fill=red!20, opacity=0.5] (8,-4) rectangle ++(-1,8);

    \node at (4.5,3.5) {Positive \\ weights};
    \node at (7.5,3.5) {Positive \\ weights};

    \draw (0,1.5*1.5) -- (0,5);
    \draw (0,4.5) -- (6,4.5);
    \draw[->] (6,4.5) -- (6,3.5);
    \draw (0,5) -- (9,5);
    \draw[->] (9,5) -- (9,0.5);
\end{tikzpicture}
    \end{subfigure}
    \hfill
    \begin{subfigure}[b]{0.49\textwidth}
    \centering 
    \begin{tikzpicture}
	
	\definecolor{darkgray176}{RGB}{176,176,176}
	\definecolor{lightgray204}{RGB}{204,204,204}
	
	\begin{axis}[
		width=0.98\textwidth,
		height=7.2cm,
		legend cell align={left},
		legend style={fill opacity=1.0, draw opacity=1, text opacity=1, draw=lightgray204, at={(0.98,0.98)}, anchor=north east},
		log basis y={10},
		tick align=outside,
		tick pos=left,
		x grid style={darkgray176},
		xlabel={Epoch},
		xmajorgrids,
		xmin=-0.3, xmax=28.3,
		xminorgrids,
		xtick style={color=black},
		y grid style={darkgray176},
		ymajorgrids,
		ymin=4.94873592530113e-05, ymax=0.0864108323693739,
		yminorgrids,
		ymode=log,
		ytick style={color=black}
		]
		
		\addplot [very thick, coltrainLoss]
		table {%
			1 0.0615461927634979
			2 0.00420070456493787
			3 0.00343826968649544
			4 0.0029564614719944
			5 0.00266809398455985
			6 0.00243896641012466
			7 0.00225792546143998
			8 0.00209647417771572
			9 0.00196406313595056
			10 0.00185166316663176
			11 0.00178271303516409
			12 0.00171307176081286
			13 0.00163076315798749
			14 0.00155334391659829
			15 0.00148399108499928
			16 0.00144764458996624
			17 0.00143185160502794
			18 0.00142006020132935
			19 0.00140999558718041
			20 0.00139509250587344
			21 0.00138975537723504
			22 0.00138598815671725
			23 0.00137908343568756
			24 0.00137096835184163
			25 0.00136371530479055
			26 0.00135861144953472
			27 0.00135608247760448
		};
		\addlegendentry{Train~\(\Loss\)}
		
		\addplot [very thick, colvalLoss]
		table {%
			1 0.0046392231701871
			2 0.00370156788036016
			3 0.00316851148938861
			4 0.00276371468199818
			5 0.00256317317827074
			6 0.00227574284516504
			7 0.00220529558462912
			8 0.00198087564823669
			9 0.00186288834043094
			10 0.00183614294420065
			11 0.00183433112373518
			12 0.00166377659102166
			13 0.00159340627930715
			14 0.00150046898406571
			15 0.00141893764441485
			16 0.00140292915066453
			17 0.00139314744160257
			18 0.00147078505445307
			19 0.00138368236799295
			20 0.00136021581982608
			21 0.00135514123898096
			22 0.00132384376077479
			23 0.00137779378217709
			24 0.00135510940262145
			25 0.0013574304444777
			26 0.00136366576918491
			27 0.00134313300067017
		};
		\addlegendentry{Val~\(\Loss\)}
		
		\addplot [very thick, colvalLossMse]
		table {%
			1 0.00360533368082422
			2 0.00302354151319193
			3 0.00262309247062015
			4 0.00227528997877476
			5 0.00202069520756567
			6 0.00185708036385344
			7 0.00187847572175156
			8 0.00167905864511654
			9 0.00158167163838124
			10 0.00146895741240229
			11 0.00142604244497536
			12 0.00137450886711954
			13 0.00136136165345594
			14 0.00109347590151734
			15 0.00111536164260546
			16 0.00107338960101964
			17 0.00114216333813305
			18 0.00121219778044841
			19 0.00116506826622416
			20 0.00111123810407769
			21 0.00104947589090134
			22 0.0010788987388165
			23 0.0012075752600913
			24 0.00113547516552991
			25 0.00111277873717011
			26 0.000979604402118876
			27 0.00114716832926739
		};
		\addlegendentry{Val~\(\Loss_{\mse}\)}
		
		\addplot [very thick, colvalLossSym]
		table {%
			1 0.000416069643887918
			2 0.000238778560610057
			3 0.000202561961555739
			4 0.000174724944934896
			5 0.000221006986018048
			6 0.000161338565877931
			7 0.000171194765637323
			8 0.000132756175761412
			9 0.000115103090254397
			10 0.000153262929012297
			11 0.000177509525500323
			12 0.000134822478862836
			13 0.000128543147587176
			14 0.000122185152001635
			15 0.000120331715264799
			16 0.000113562399018999
			17 0.000118327720638284
			18 0.000151414574776106
			19 0.000121458247237739
			20 0.000116553216116216
			21 0.000114712204740327
			22 0.000103599205513371
			23 0.000100738286456399
			24 0.000123695838883835
			25 0.000105624517059035
			26 0.000107118797204662
			27 0.000103961901464921
		};
		\addlegendentry{Val~\(\lambda_{\sym} \, \Loss_{\sym}\)}
		
		\addplot [very thick, colvalLossIneq]
		table {%
			1 0.000617819855228877
			2 0.000439247807982834
			3 0.000342857055823109
			4 0.00031369975428399
			5 0.00032147098551163
			6 0.000257323912499801
			7 0.000155625099120485
			8 0.000169060827452048
			9 0.000166113610868656
			10 0.000213922604516956
			11 0.00023077915402631
			12 0.000154445244967262
			13 0.000103501476729544
			14 0.000284807929796525
			15 0.00018324428763224
			16 0.000215977147935848
			17 0.000132656383033854
			18 0.000107172698378606
			19 9.71558546129271e-05
			20 0.000132424502585126
			21 0.000190953143213735
			22 0.000141345815213308
			23 6.94802344841568e-05
			24 9.59383989093155e-05
			25 0.00013902718838751
			26 0.000276942570163162
			27 9.20027681025918e-05
		};
		\addlegendentry{Val~\(\lambda_{\ineq} \, \Loss_{\ineq}\)}
	\end{axis}
	
\end{tikzpicture}
    \end{subfigure}
    \caption{
    Left: Network architecture for Kohn--Strang--Dolzmann example. 
    Right: Learning curves based on the loss function \(\Loss\) from \eqref{eq:loss} for a single network initialisation. 
    For the validation loss, the individual contributions \(\Loss_{\mse}\), \(\lambda_{\sym} \, \Loss_{\sym}\) and \(\lambda_{\ineq} \, \Loss_{\ineq}\) to the overall loss \eqref{eq:loss} are depicted. 
    } \label{fig:Architecture_Loss_KSD}
\end{figure}
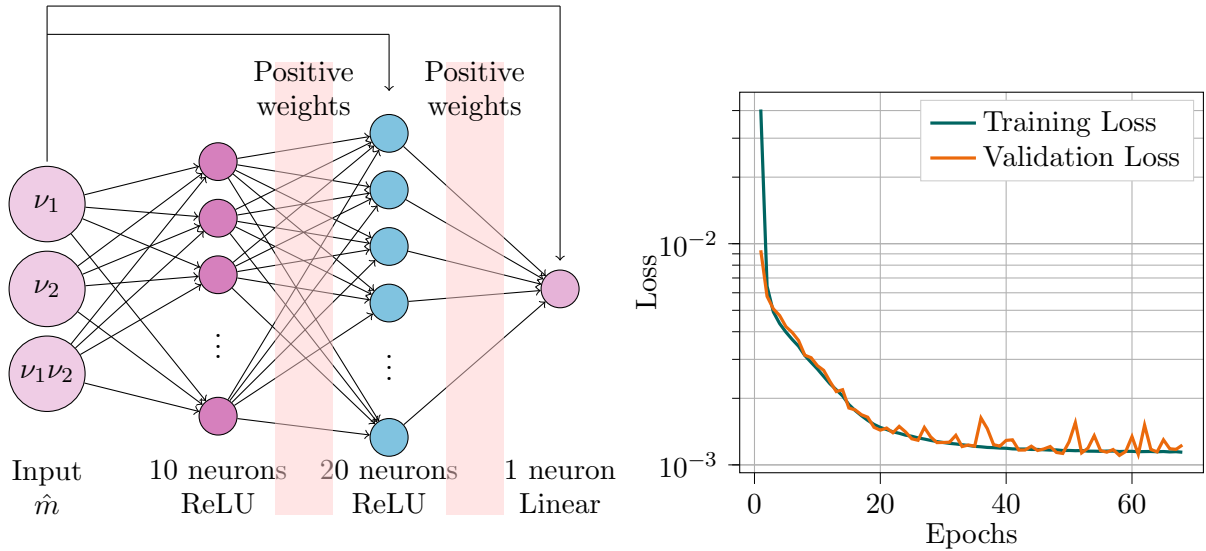

On average, for a single realisation, the training process requires \(6 \pm  1\) minutes to complete \(32 \pm 7\) epochs. 
The final training loss is \(1.3 \times 10^{-3}\pm 2.1  \times 10^{-4}\),  while the validation loss reaches \(1.3  \times 10^{-3} \pm 2.3 \times 10^{-4}\). 
The learning curves for both training and validation losses on one exemplary realisation are presented in \cref{fig:Architecture_Loss_KSD}.
In addition to the loss function \eqref{eq:loss}, the individual contributions to the overall validation loss \(\Loss\), i.e.~\(\Loss_{\mse}\) \eqref{eq:Lmse}, \(\Loss_{\ineq}\) \eqref{eq:lossineq} and \(\Loss_{\sym}\) \eqref{eq:losssym}, including the scaling by the penalisation parameters, are plotted, showing that the symmetry as well as the inequality property is taken into account and improved during the training process. In particular, the contributions of \(\Loss_{\ineq}\) and \(\Loss_{\sym}\) are significantly smaller than that of \(\Loss_{\mse}\), substantiating the weak enforcement of the corresponding properties.

\begin{table}[h]
\centering
\begin{tabular}{@{}c|c|c@{}}
Mean err & Rel quad err & Rel max err \\ \midrule
\(0.022 \pm 1.8 \times 10^{-3}\)  & \(0.021 \pm 1.6 \times 10^{-3}\) & \(0.030 \pm 5.3 \times 10^{-3}\) \\ 
\end{tabular}
\caption{Average of prediction errors over twenty network realisations.}
\label{tab:Table_KSD}
\end{table}

\begin{figure}[htbp]
	\centering
	\input{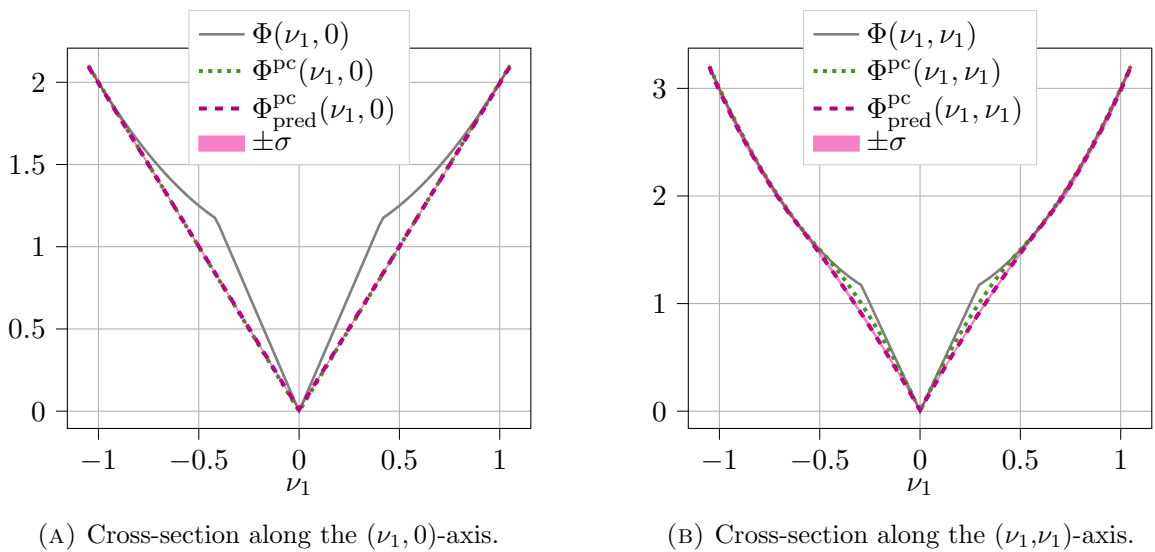}
	\caption{
		Comparison of the predicted polyconvex envelope \(\Phi^{\pc}_\pred\) (averaged over twenty network realisations) and the analytical polyconvex envelope \(\Phi^{\pc}\) for the Kohn--Strang--Dolzmann function along two cross-sections.
		The standard deviation of the twenty predictions is marked by \(\sigma\).
	}
	\label{fig:KSD_1D}
\end{figure}

\Cref{tab:Table_KSD} presents the approximation errors of the predictions, where the analytical polyconvex envelope is used as a ground truth. 
These errors are computed on a uniform $100 \times 100$ discretisation of the domain $[-\overline{\nu},\overline{\nu}]^2$. 
The results indicate that the predicted envelopes deviate by about 2\% to 3\% from the analytical polyconvex envelope \(\Phi^{\pc}\), demonstrating an accuracy sufficient for the intended engineering applications. The approximation quality is further illustrated in \cref{fig:KSD_1D}, where two one-dimensional cross-sections of the predicted polyconvex envelope are depicted. It is important to note that the neural network successfully captures the kink of the function \(\Phi^{\pc}\) at the point $(0,0)$ with high accuracy.
Further, it is notable that the neural network based approximation captures the nonconvexity along the diagonal cross-section of the envelope, rendering it a consistent polyconvex function.
It is observable that the standard deviation \(\sigma\) of the different network realisations is quite negligible when it comes to the approximation accuracy. 
Specifically, the neural network implemented in this example contains \num{344} parameters, comprising both weights and biases. 
With the same number of parameters using the conventional approach, it would only be possible to store \(19 \times 19\) grid values, and then approximate the value at any point by interpolation. 
While this standard method remains feasible for simple cases, it can quickly become intractable for more complex scenarios, such as parameter-dependent functions or real engineering problems.

\begin{rem}
In addition to the values of the polyconvex envelopes, most engineering applications also require their derivatives.
We briefly point out two straight forward approaches for predicting these derivatives. 
The first approach involves differentiating the predictions of the neural network directly using standard finite difference methods. 
However, preliminary results indicate that this method suffers from significant drawbacks owing to the lack of continuity in the network evaluations. An alternative is to train a separate neural network to learn the derivatives of the polyconvex envelopes. Initial numerical experiments on this example demonstrate that this can be achieved with a relatively simple network architecture, consisting of two hidden layers with 10 and 20 neurons, ReLU activation functions and a two-dimensional output with linear activation functions. 
However, we do not pursue this direction further as it is beyond the scope of this paper.
\end{rem}

\subsection{Saint~Venant--Kirchhoff Energy Density}
Building on the previous two-dimensional example, we now move to a function well-known in elasticity simulations and illustrate the feasibility of our approach in three spatial dimensions.
Specifically, we consider the Saint~Venant--Kirchhoff model in three spatial dimensions with determinant constraints given by
\begin{equation} \label{eq:WSTVKdetconstr}
    W(F) = 
    \begin{cases}
    	\frac{\mu}{4} \, \lvert F^T F - \mathbb{I}\rvert^2 + \frac{\lambda}{8} \left(\lvert F\rvert^2 - 3 \right)^2 & \text{ if } \det F > 0, \\
    	\infty & \text{ else.}
    \end{cases}
\end{equation}
In contrast to the unconstrained case \eqref{eq:STVK}, where the polyconvex envelope is known analytically \cite{LeRao:1995:qes}, no analytical polyconvex envelope is available in the determinant-constrained setting, making numerical approximations necessary.
Accordingly, its evaluation within an offline--online framework necessitates the storage of a large set of values on a three-dimensional grid, followed by interpolation. 
This example thus allows to illustrate a practical application of our approach, namely the efficient compression of a precomputed polyconvex envelope by a properties-preserving neural network, which in turn enables accurate and computationally efficient online prediction.

In order to learn \(\Phi^{\pc}\), we focus on the signed singular value reformulation of \eqref{eq:WSTVKdetconstr}, which is given by
\begin{equation} \label{eq:PhiSTVKdetconstr}
	\Phi(\hat{\nu}) = 
	\begin{cases}
		\frac{\mu}{4} \, \sum_{i = 1}^{3} (\nu_i^2 - 1)^2 + \frac{\lambda}{8} \, (\lvert\hat{\nu}\rvert^2 - 3)^2 & \text{ if } \nu_1 \,\nu_2\,\nu_3 > 0, \\
		\infty & \text{ else.}
	\end{cases}
\end{equation}
Within this formulation, we consider the material parameters \(\mu = 0.4\) and \(\lambda = 0.4\) for the numerical experiments. 

\subsubsection{Network Architecture and Learning Data}
We implement a FICNN as described in \cref{fig:FICNN}, noting that the input layer is now of dimension \(k_3 = 7\), with three hidden layers consisting of 20, 25 and 30 neurons, respectively, leading to a total of \num{1888} trainable parameters. 
To generate the learning data, we approximate the polyconvex envelope of the three-dimensional function \(\Phi\) \eqref{eq:PhiSTVKdetconstr} by the SVPC~LP algorithm described in \cref{sec:SVPC}, with the parameters \(\delta = 0.047\) and \(r = 1.5\). 
The learning domain for \(\hat{\nu}\) is defined as the box \([\nu_{\min}, \nu_{\max}]^3\) with \(\nu_{\min} = 0.4\) and \(\nu_{\max} = 1.4\). 
The learning data is generated by uniformly discretising each coordinate axis within the box by \(50\) points and considering all possible \(\Pid\) permutations of this discretisation, resulting in \num{1000000} data points covering the entire signed singular value space, half of which yield finite function values of \(\Phi\). 
The penalty parameters in the loss function \eqref{eq:loss} are chosen as \(\lambda_{\sym} = 1.0\) and \(\lambda_{\ineq} = 2.0\).

\subsubsection{Numerical Results}
On average, for a single realisation, the training process requires \(17 \pm 6\) minutes to complete \(32 \pm 10\) epochs. 
The final training loss is \(1.15 \times 10^{-5} \pm 2.4 \times 10^{-6}\), while the final validation loss reaches \(1.73 \times 10^{-5} \pm 1.2 \times 10^{-5}\). 
The mean squared error \eqref{eq:Lmse} between the neural network approximation and the learning data, i.e.~reference polyconvex envelope \(\Phi^{\pc}_{\delta}\), accounts to \(5.03 \times 10^{-6} \pm 1.3 \times 10^{-6}\). 
The mean error on the learning data is \(0.002 \pm 2.3 \times 10^{-4}\), the relative quadratic error is \(0.020 \pm 2.5 \times 10^{-3}\) and the relative maximum error is \(0.027 \pm 5.9 \times 10^{-3}\).
\begin{figure}
	\centering
	\input{Pictures/STVK/STVK_slices_3d.tex}
	\caption{
		Cross-sections for the three-dimensional Saint~Venant--Kirchhoff model \eqref{eq:PhiSTVKdetconstr}.
		\(\Phi^{\pc}_{\pred}\) is averaged over ten network realisations.
		The reference polyconvex envelope \(\Phi^{\pc}_{\delta}\) is generated by the SVPC~LP algorithm, see \cite{NPPW2024} and \cref{sec:SVPC}. 
		The evaluation of \(\nu_1 \in [0.4, 1.4]\) involves \num{200} points on each cross-section that do not belong to the learning dataset.
		Note the different scalings of the ordinate axes in the three plots due to the different axial cross-sections.
	}
	\label{fig:STVK3Dcrosssections}
\end{figure}

\cref{fig:STVK3Dcrosssections} shows three cross-sections, \((\nu_1, \nu_1, \nu_1)\)-axis (left), \((\nu_1, \nu_1, 1)\)-axis (middle) and \((\nu_1, 1, 1)\)-axis (right), for the predictions of the neural network model, corresponding to triaxial, biaxial and uniaxial deformations in the signed singular value space, respectively. 
The reference polyconvex envelope \(\Phi^{\pc}_{\delta}\) in \cref{fig:STVK3Dcrosssections} is computed by the SVPC~LP algorithm, see \cite{NPPW2024} and illustrated in \cref{sec:SVPC} using the discretisation parameters \(r = 1.5\) and \(\delta = 0.025\), and evaluated for \(\nu_1 \in [\nu_{\min}, \nu_{\max}] = [0.4, 1.4]\) uniformly discretised in \num{200} points. 
The resulting points on the cross-sections are not part of the learning data set, thus demonstrating the interpolation quality of the network.
On all of the three illustrated cross-sections, the maximum absolute error is found to be equal to \(0.0092\), which corresponds to a relative maximum error of \(2\%\). 
This error is of the same order as the error observed in the learning data for the entire box \([\nu_{\min}, \nu_{\max}]^3\). 
These numerical results underline the compression potential of the neural-network-based representation approach. 
In particular, the proposed neural network compresses the large learning data set, consisting of \(\sim 100^3\) point values of \(\Phi^{\pc}_{\delta}\), into only \num{1888} trainable parameters, corresponding to merely a \(12^3\) grid.

\subsection{A Parametric Example: The Generalised Kohn--Strang--Dolzmann Function}
\label{sec:GKSD}
After the previous initial benchmark examples, we consider a more complex case, i.e.~a two-parameter-dependent family of functions. 
The following example was studied in \cite{AllFranc,Zhang} and modified to achieve continuity. 
We consider the function $W \colon \R^{2\times 2} \times \R_{+} \times \R_{+}\setminus\{0\} \rightarrow \R$, defined as
\begin{equation*}
    W(F; \lambda, \alpha)=
    \begin{cases}
        \lambda + \alpha \, \lvert F \rvert^2 & \text{if } \, \lvert F \rvert \geq \sqrt{\frac{\lambda}{\alpha}} \, (\sqrt{2}-1), \\ 
        2\sqrt{2\, \lambda \,\alpha} \, \lvert F \rvert & \text{else},
    \end{cases}
\end{equation*}
where $\lvert F \rvert \coloneqq (\sum_{i,j=1}^{d} F_{ij}^2)^{1/2}$. The polyconvex envelope of $W$ is known analytically in closed form and reads
\begin{equation*}
    W^{\pc}(F; \lambda, \alpha)=
    \begin{cases}
        \lambda + \alpha \, \lvert F \rvert^2 & \text{if } \, \rho(F) \geq \sqrt{\frac{\lambda}{\alpha}}, \\ 
        2 \sqrt{\lambda \, \alpha} \, \rho(F) - 2 \alpha \, \lvert {\det}(F) \rvert & \text{else},
    \end{cases}
\end{equation*}
where $\rho(F) \coloneqq \sqrt{\lvert F \rvert^2 + 2\, \lvert {\det}(F) \rvert}$. 
The functions $W$ and $W^{\pc}$ are isotropic, and, rewritten in terms of the signed singular values, they reduce to $\Phi$, $\Phi^{\pc}\colon \R^2 \times \R_{+} \times \R_{+} \setminus\{0\}\to \R$ with 
\begin{equation}
    \Phi(\hat{\nu}; \lambda, \alpha) =
    \begin{cases}
        \lambda + \alpha \, (\nu_1^2 +\nu_2^2) & \text{if } \, \sqrt{\nu_1^2 +\nu_2^2} \geq \sqrt{\frac{\lambda}{\alpha}} \, (\sqrt{2}-1),\\
        2\sqrt{2\, \lambda \, \alpha} \, \sqrt{\nu_1^2 +\nu_2^2} & \text{else},
    \end{cases}
\end{equation}
and 
\begin{equation}
    \Phi^{\pc}(\hat{\nu}; \lambda, \alpha)=
    \begin{cases}
        \lambda + \alpha(\nu_1^2 +\nu_2^2) & \text{if } \ \rho(\hat{\nu}) \geq \sqrt{\frac{\lambda}{\alpha}},\\ 
        2 \sqrt{\lambda \, \alpha} \, (\lvert \nu_1 \rvert + \lvert \nu_2 \rvert) - 2 \alpha \, \lvert \nu_1 \nu_2 \rvert & \text{else},
    \end{cases}
\label{eq:PhipcGKSD}
\end{equation}
where $\rho(\hat{\nu}) = \lvert \nu_1 \rvert + \lvert \nu_2 \rvert$.

\subsubsection{Network Architecture and Learning Data}
For this parameter-dependent example, we implement a PICNN where each path consists of three hidden layers with 10, 20 and 20 neurons, respectively. 
The convex input, denoted by $\hat{m}$, represents the minors of the signed singular values, i.e.~$\m(\hat{\nu})$, while the nonconvex inputs, denoted as vector $\zeta$, correspond to the two parameters $\lambda$ and $\alpha$. 
The learning domain for the parameters \(\lambda\) and \(\alpha\) is \([1,2]\), while for the signed singular value input \(\hat{\nu}\) the learning domain is defined as \([-\overline{\nu}, \overline{\nu}]^2\), with \(\overline{\nu}=1.5\), ensuring in particular that \(\overline{\nu} \geq \sqrt{\lambda/\alpha}\).
For the training dataset, the set of points in the signed singular value space is obtained by discretising each axis into \num{251} points, so that the point \(0\) is included, with a local refinement towards the origin using a quadratic transformation, as in the example in \cref{sec:KSD}.
The parameter $\lambda$ takes values from the discrete set $\{1,1.2, \ldots, 1.8,2\}$ and $\alpha$ follows the same parameter discretisation.
These discretisations lead to a total of \num{2268036} points in the learning dataset. 
For the validation dataset, $30 \%$ of these points are randomly selected, leaving $70 \%$ for the training data set. 
The target values are computed using the analytical function \(\Phi^{\pc}\), cf.~\eqref{eq:PhipcGKSD}. 

\begin{figure}[htbp]
	\centering
	\begin{subfigure}[b]{0.48\textwidth}
		\centering 
		\input{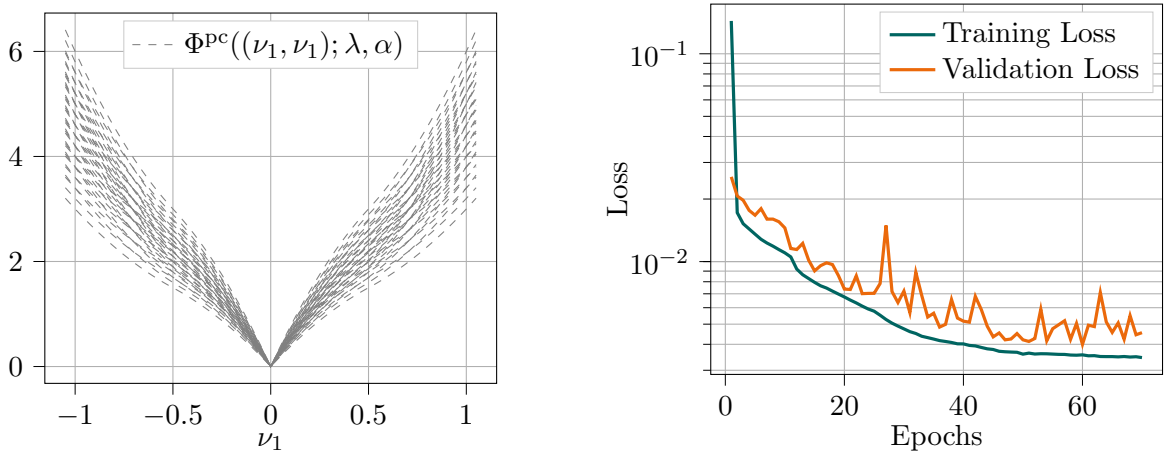}
	\end{subfigure}
	\hfill
	\begin{subfigure}[b]{0.48\textwidth}
		\centering 
		\begin{tikzpicture}
	
	\definecolor{darkgray176}{RGB}{176,176,176}
	\definecolor{darkorange}{RGB}{255,140,0}
	\definecolor{lightgray204}{RGB}{204,204,204}
	\definecolor{royalblue}{RGB}{65,105,225}		
	\begin{axis}[
		width=0.98\textwidth,
		legend cell align={left},
		legend style={
			fill opacity=1,
			draw opacity=1,
			text opacity=1,
			at={(0.98, 1.15)},
			anchor=north east,
			draw=lightgray204
		},
		log basis y={10},
		tick align=outside,
		tick pos=left,
		x grid style={darkgray176},
		xlabel={Epoch},
		xmajorgrids,
		xminorgrids,
		xtick style={color=black},
		y grid style={darkgray176},
		ymajorgrids,
		yminorgrids,
		ymode=log,
		ytick style={color=black},
		xmin=0, xmax=22,
		ymin=0.000183466704535238, ymax=0.211860984124631,
		]
		
		\addplot [very thick, coltrainLoss]
		table {%
			1 0.153761187739518
			2 0.0303913935082953
			3 0.0228834660537593
			4 0.0204066155387318
			5 0.018630867126508
			6 0.0170065793027544
			7 0.0157749658670185
			8 0.0149398392433314
			9 0.0139816728920236
			10 0.0127141854870025
			11 0.0119877595935179
			12 0.0116951202067331
			13 0.0109932458597395
			14 0.0103574499787148
			15 0.0099440672209452
			16 0.00959053828291378
			17 0.00939047616343183
			18 0.00926407207113255
			19 0.00903290485369969
			20 0.00895665810952251
			21 0.00884121939836909
		};
		\addlegendentry{Train~\(\Loss\)}
		
		\addplot [very thick, colvalLoss]
		table {%
			1 0.036585188949392
			2 0.025548032407523
			3 0.0196878633781048
			4 0.0176211694781269
			5 0.0166462003222453
			6 0.020103440466591
			7 0.0140699427398329
			8 0.0132878561632042
			9 0.0119712050492828
			10 0.0119463860006855
			11 0.0111966119673923
			12 0.0112219635860222
			13 0.00907051696697165
			14 0.00988392917096684
			15 0.0100658086910879
			16 0.00822720611248955
			17 0.00823966557024047
			18 0.0132493553441301
			19 0.00911998441888107
			20 0.00905952648187031
			21 0.00962908531607539
		};
		\addlegendentry{Val~\(\Loss\)}
		
		\addplot [very thick, colvalLossMse]
		table {%
			1 0.0279059029215722
			2 0.0201716209910843
			3 0.0140486565975671
			4 0.0109331522505152
			5 0.0105121959916941
			6 0.0166212955296676
			7 0.00758133680203528
			8 0.00854876640976806
			9 0.00725682275304747
			10 0.00639400010334275
			11 0.00777956851486513
			12 0.00726799270922941
			13 0.00564951962888826
			14 0.00634597047344453
			15 0.00542872233652033
			16 0.00511278063969357
			17 0.00596114688709599
			18 0.0101060459728231
			19 0.00688541018081736
			20 0.00595998086072455
			21 0.00680705685066203
		};
		\addlegendentry{Val~\(\Loss_{\mse}\)}
		
		\addplot [very thick, colvalLossSym]
		table {%
			1 0.00510895829242864
			2 0.00424201238936634
			3 0.00361433193915668
			4 0.00378092671369021
			5 0.0039917446512769
			6 0.00312252053342216
			7 0.00252986769850896
			8 0.00330505442587995
			9 0.00291350785802799
			10 0.00277596629074311
			11 0.00283892199526891
			12 0.0031918521261479
			13 0.00187432179832064
			14 0.00286811484863995
			15 0.00312834673684184
			16 0.00181838205421335
			17 0.0016694473561034
			18 0.0027391326053216
			19 0.00185823773076698
			20 0.00258861439915751
			21 0.00256923751543292
		};
		\addlegendentry{Val~\(\lambda_{\sym} \, \Loss_{\sym}\)}

		\addplot [very thick, colvalLossIneq]
		table {%
			1 0.00357032778453877
			2 0.00113439901934136
			3 0.00202487483941007
			4 0.00290709050334615
			5 0.00214225968713764
			6 0.000359624395905952
			7 0.00395873824537701
			8 0.00143403533078188
			9 0.00180087442706439
			10 0.00277641960568749
			11 0.000578121451340066
			12 0.000762118753579543
			13 0.00154667555112682
			14 0.000669843855192465
			15 0.00150873961081498
			16 0.00129604341276639
			17 0.000609071329658956
			18 0.000404176774158872
			19 0.000376336504877014
			20 0.000510931217214589
			21 0.000252790949057872
		};
		\addlegendentry{Val~\(\lambda_{\ineq} \, \Loss_{\ineq}\)}
	\end{axis}
	
\end{tikzpicture}
	\end{subfigure}
	\caption{
		Left: Training data for the generalised Kohn--Strang--Dolzmann example.
		One-dimensional cross-section along the diagonal \((\nu_{1}, \nu_{1})\)-axis of the analytical polyconvex envelopes.
		Envelopes are evaluated at the training points in parameter space, i.e.~for \((\lambda, \alpha) \in \{1, 1.2, \dotsc, 1.8, 2 \}^2\). 
		Right:
		Learning curves based on the loss function \(\Loss\) from \eqref{eq:loss} for a single network initialisation. 
	}
	\label{fig:Training_Loss_GKSD}
\end{figure}

\Cref{fig:Training_Loss_GKSD} illustrates a one-dimensional representation of the target values for the training dataset along the diagonal $(\nu_1, \nu_1)$-axis. 
For the loss function \(\Loss\) from \eqref{eq:loss}, we choose the penalty parameters $\lambda_{\ineq}=50$ and $\lambda_{\sym}=20$.

\subsubsection{Numerical Results}
On average, for a single realisation, the training process requires \(27 \pm 5\) minutes to complete \(22 \pm 4\) epochs. 
The final training loss is \(8.07 \times 10^{-3} \pm 1.6 \times 10^{-3}\), while the validation loss reaches \(8.57 \times 10^{-3} \pm 1.4 \times 10^{-3}\). 
The learning curves for one selected realisation are plotted in \cref{fig:Training_Loss_GKSD}.

\begin{figure}[htbp]
	\centering
	\begin{tikzpicture}

	\definecolor{darkgray176}{RGB}{176,176,176}
	
	\pgfmathsetmacro{\xshift}{0.5} 
	
 	\pgfplotsset{
		standard axis/.style={
			width=0.27\textwidth,
			axis equal image,
			scale only axis,
			colorbar horizontal, 
			colorbar style={
				at={(0.5,-0.3)},
				anchor=north,
				height=0.2cm,
			},
			colorbar style={
				xticklabel style={
					/pgf/number format/fixed,
					/pgf/number format/precision=3,
					/pgf/number format/fixed zerofill
				},
				scaled ticks=false,
			},	
			colormap name=viridis,
			grid=both, 
			grid style={dashed, gray!30},
			scatter,
			only marks,
			mark size=2.0pt,
			tick align=center,
			xtick={1, 1.2, ..., 2.0},
			ytick={1, 1.2, ..., 2.0},
			xticklabel style={/pgf/number format/fixed, /pgf/number format/precision=1, /pgf/number format/fixed zerofill},
			yticklabel style={/pgf/number format/fixed, /pgf/number format/precision=1, /pgf/number format/fixed zerofill},	
		}
	}
	
	\begin{axis}[
		name=plot1,
		standard axis,
		xlabel={\(\lambda\)}, 
		ylabel={\(\alpha\)},
		ylabel style={yshift=-7pt},		title={Mean err},
		]
		\addplot[scatter, only marks, scatter src=explicit, mark size=2.0pt] 
		table[header=false, x index=0, y index=1, meta index=2] 
		{Pictures/GenKSD/PhipcNNGKSD_paramspace_meanabs_errors.dat};

	\end{axis}
	
	\begin{axis}[
		name=plot2,
		standard axis,
		at={(plot1.right of north east)},
		anchor=left of north west,
		xshift=\xshift cm,
		xlabel={\(\lambda\)},
		yticklabels={},
		title={Rel quad err}
		]
		\addplot[scatter, only marks, scatter src=explicit, mark size=2.0pt] 
		table[header=false, x index=0, y index=1, meta index=2]
		{Pictures/GenKSD/PhipcNNGKSD_paramspace_relquad_errors.dat};
	\end{axis}
	
	\begin{axis}[
		name=plot3,
		standard axis,
		at={(plot2.right of north east)},
		anchor=left of north west,
		xshift=\xshift cm,
		xlabel={\(\lambda\)},
		yticklabels={},
		title={Rel max err}
		]
		\addplot[scatter, only marks, scatter src=explicit, mark size=2.0pt] 
		table[header=false, x index=0, y index=1, meta index=2] 
		{Pictures/GenKSD/PhipcNNGKSD_paramspace_relmax_errors.dat};
	\end{axis}
	
\end{tikzpicture}
	\caption{
		Errors over parameter space for the generalised Kohn--Strang--Dolzmann example. 
		Left: Mean errors.
		Middle: Relative quadratic errors.
		Right: Relative max errors.
		All errors are computed for an average over ten network realisations with respect to the analytical polyconvex envelopes for different values $\lambda$ and $\alpha$.}
	\label{fig:Error_GKSD}
\end{figure}
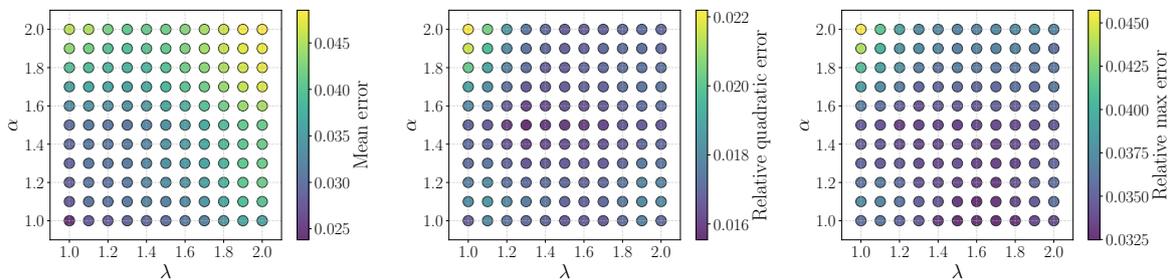

\begin{figure}[htbp]
	\centering
	\input{Pictures/GenKSD/Phi_GKSD_1d_1.7_1.3_ReLU.tex}
	\caption{
		Comparison along two cross-sections of the predicted polyconvex envelope \(\Phi^{\pc}_\pred\) (averaged over ten network realisations) and the analytical polyconvex envelope \(\Phi^{\pc}\) for the Generalised Kohn--Strang--Dolzmann function with $\lambda=1.7$ and $\alpha=1.3$.
		The standard deviation of the ten predictions is marked by $\sigma$.}
	\label{fig:GKSD_1D_2}
\end{figure}

\begin{figure}[htbp]
	\centering
	\input{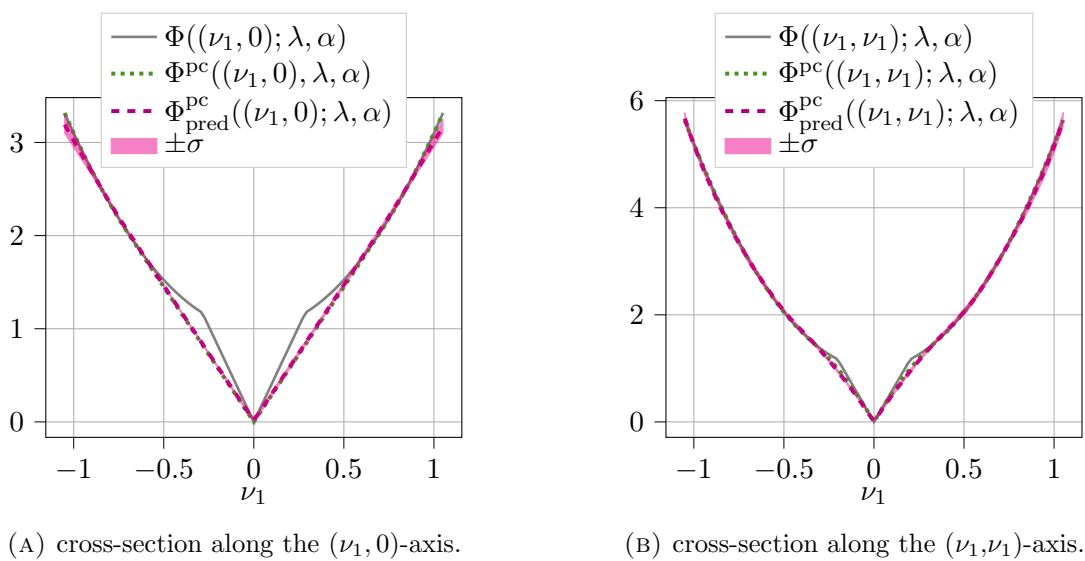}
    \caption{
    	Comparison along two cross-sections of the predicted polyconvex envelope \(\Phi^{\pc}_\pred\) (averaged over ten network realisations) and the analytical polyconvex envelope \(\Phi^{\pc}\) for the Generalised Kohn--Strang--Dolzmann function with $\lambda=2.3$ and $\alpha=2.4$; these parameters are located outside the parameter learning domain \([1,2]\). 
		The standard deviation of the ten predictions is marked by $\sigma$.}
	\label{fig:GKSD_1D_3}
\end{figure}

\Cref{fig:Error_GKSD} illustrates the relative errors between the predictions and the analytical polyconvex envelopes for different sets of parameters $(\lambda, \alpha)$, including both values present in the training set and those outside of it. These errors are computed on a uniform $(100 \times 100)$-discretisation of the domain $[-\overline{\nu},\overline{\nu}]^2$.
Across all considered parameter sets, the errors remain consistently low, ranging between 2\% and 4\%.  Additionally, \cref{fig:GKSD_1D_2} provides one-dimensional cross-sections of these predictions. Notably, the numerical experiments indicate that errors are primarily localised at the domain boundaries, while the neural network successfully captures the kink at the point $(0,0)$ with an accuracy sufficient for typical engineering applications. 
Furthermore, \cref{fig:GKSD_1D_2} demonstrates that the neural network is capable of making accurate predictions even for values of $\hat{\nu}$ outside the training domain, defined as $[-1.5, 1.5]^2$. In particular, the predictions remain accurate for $\hat{\nu}$ in $[-1.8, 1.8]^2$. 
This highlights the network's ability to extrapolate reliably in the signed singular values argument. 
The behaviour in the extrapolation regime is particularly relevant in finite element simulations, where unreasonably large strains may arise during Newton iterations.

A further extrapolation analysis shows that the implemented neural network performs well even for parameter pairs $\{\lambda, \alpha\}$ in the range $[1, 2]^2 \cup [1.5, 2.5]^2$, i.e.~also for pairs outside the training domain, thereby indicating a degree of extrapolation capability in terms of parameters as well (but inevitably with a loss in the precision). An example of such behaviour is illustrated in \cref{fig:GKSD_1D_3}. Beyond this range, however, the prediction quality deteriorates noticeably. Nevertheless, extrapolation in these parameters is less critical than extrapolation in the signed singular values since the parameters are either known (and possibly constant) or bounded (see e.g.~\cref{sec:numericalexperimentsisodamage}) for a given model. 

Apart from the extrapolation capacity, the crucial aspect here is its compression capability within the training range, which will become even more apparent in higher-dimensional parameter spaces.
In particular, the neural network implemented in this example contains \num{3291} parameters, comprising both weights and biases. 
With the same number of parameters using the conventional approach, it would only be possible to store a $19 \times 19$ grid, with 3 values for each $\lambda$ and $\alpha$.  With this limited grid size, it is evident that the full range of data cannot be recovered with the same precision as is achievable with the neural networks presented here.

\section{Engineering Example: Isotropic Damage Problem} \label{sec:numericalexperimentsisodamage}
Having validated our approach with (mathematical) benchmark problems, we return to the isotropic damage problem introduced in \cref{sec:damage}.
In order to fit this model into our setting, we rephrase it in signed singular value formulation.  

\subsection{Reformulation of the Isotropic Damage Problem} \label{sec:isodamage}
The function \(W\) from \eqref{eq:Wdamage} is dependent on \(d\times d + 1\) parameters and on the \(d \times d\)-deformation gradient.
Due to the isotropy of \(W\), it can be recast into a signed singular value formulation utilising \(\hat{\nu} \in \R^d\). 
To stress the dependence on the signed singular values, the function \(\varphi\) is now employed to rewrite \eqref{eq:psi1minusD} as 
\begin{equation*}
	\varphi(\hat{\nu}, \alpha) = (1 - D(\alpha)) \, \varphi^{0}(\hat{\nu}).
\end{equation*}
Within this formulation, the function \(\varphi^{0}\) denotes the signed singular value formulation of the isotropic undamaged energy density \(\psi^{0}\) from \eqref{eq:psi1minusD}, i.e.~\(\varphi^{0}\) and \(\psi^{0}\) are related by 
\begin{equation*}
	\varphi^0(\hat{\nu}) = \psi^0(\diag(\hat{\nu})) \qquad\quad \text{and} \qquad\quad \psi^0(F) = \varphi^0(\nu(F)).
\end{equation*}
Consequently, the pseudo-time incremental energy density \(W\) in the signed singular value formulation, denoted by \(\Phi\), for the time step \( k + 1 \) reads
\begin{align} \label{eq:Phidamage}
	\begin{aligned}
		\Phi(\hat{\nu}_{k+1}; \hat{\nu}_{k}, \alpha_{k}) 
		& = \varphi(\hat{\nu}_{k+1},p(\hat{\nu}_{k + 1}; \alpha_{k})) - \varphi(\hat{\nu}_k,\alpha_k) \\
		&\qquad + p(\hat{\nu}_{k + 1}; \alpha_{k})\, D(p(\hat{\nu}_{k + 1}; \alpha_{k})) - \alpha_k \, D(\alpha_k) 
		- \overline{D}(p(\hat{\nu}_{k + 1}; \alpha_{k})) + \overline{D}(\alpha_k).
	\end{aligned}
\end{align}
The parameter dependence in this function can be given the following interpretation. 
The scalar parameter \(\alpha_{k}\) plays still the same role as the internal variable, while the vector \(\hat{\nu}_{k} \in \R^{d}\) belongs to the signed singular values of the deformation gradient \(F_{k}\) from the previous pseudo-time step.
Within \eqref{eq:Phidamage}, the internal variable evolution is written explicitly using the path function in the signed singular value formulation as 
\begin{equation} \label{eq:pathfunctionIsotropic}
	\alpha_{k+1} 
    =
    p(\hat{\nu}_{k + 1}; \alpha_{k}) = 
	\begin{cases}
		\varphi^{0}(\hat{\nu}_{k + 1}) & \text{if } \, \varphi^{0}(\hat{\nu}_{k + 1}) > \alpha_{k}, \\
		\alpha_{k} & \text{else.}
	\end{cases}
\end{equation}
The formulation as stated in \eqref{eq:Phidamage} significantly reduces dimensionality in both the \(\hat{\nu}_{k+1}\) argument as well as the \(\hat{\nu}_{k}\) parameter dependence, opening the possibility for efficient parameter-dependent polyconvexification of the function \(\Phi(\hat{\nu}_{k+1}; \hat{\nu}_{k}, \alpha_{k})\) in the argument \(\hat{\nu}_{k+1}\). 
For the neural network, this leads to a reduction in parameter space from \(d \times d + 1\) to \(d + 1\) dimensions, and a reduction in the minors input argument from dimension \(K_d\) to \(k_d\), corresponding to the dimension of the vector \(\m(\hat{\nu}_{k+1})\). 

The damage parameters $d_0$ and $d_\infty$ are set to \(d_{0}=0.5, d_{\infty}=0.99\) in our numerical experiments. The Lamé constants $\lambda$ and $\mu$ of the materials in \eqref{eq:STVK} and \eqref{eq:NH} are set to $\lambda=0$, $\mu=0.5$, as in \cite{KNMPPB2022}.

\begin{rem}
Let us consider $\alpha_{\infty}$ and $k_0$ such that for all $k \geq k_0$ it holds $\alpha_k \geq \alpha_{\infty}$ and by evolution it holds $\alpha_{k+1} \geq \alpha_k$. 
For the choice of damage function $D$  in \eqref{eq:damagefunction}, we have for $\alpha_{\infty}$ large enough that $D(\alpha_k) \approx d_\infty$ and $\overline{D}(\alpha_k) \approx \alpha_k \, d_\infty$ for \(d_{\infty}\) marking the asymptotic damage limit. 
In such a case, \eqref{eq:Phidamage} can be rewritten as
\begin{align*}
	\begin{aligned}
		\Phi(\hat{\nu}_{k+1}; \hat{\nu}_{k}, \alpha_{k}) 
		& \approx (1-d_\infty) \, \varphi^0(\hat{\nu}_{k+1}) - 
        (1-d_\infty) \, \varphi^0(\hat{\nu}_{k})\\
		&\qquad\qquad\qquad\qquad\qquad\qquad + \alpha_{k + 1} \, d_\infty - \alpha_k \, d_\infty - \alpha_{k+1} \, d_\infty + \alpha_k \, d_\infty \\
	    & \approx (1-d_\infty) \, (\varphi^0(\hat{\nu}_{k+1})- \varphi^0(\hat{\nu}_{k})).
    \end{aligned}
\end{align*}
Consequently, for $\alpha_k \geq \alpha_{\infty}$, the energy becomes independent of $\alpha_k$. 
This observation allows to train the neural network for the parameter $\alpha_{k} \in [0, \alpha_\infty]$ and to consider $\alpha_{k} = \alpha_\infty$ for the predictions in the case $\alpha_{k} \geq \alpha_\infty$, which drastically reduces the computational effort.
For the choice of parameters $d_0$ and $d_\infty$, we choose $\alpha_\infty=4$ in the numerical experiments noting that for $\alpha_k \geq \alpha_\infty$, we can perform the estimate 
\begin{equation*}
    \lvert D(\alpha_\infty) - D(\alpha_k) \rvert \leq \left\lvert \exp\left(-\frac{\alpha_\infty}{d_0}\right) - \exp\left(-\frac{\alpha_k}{d_0}\right) \right\rvert \leq \exp\left(-\frac{\alpha_\infty}{d_0}\right) = \exp(-8) \approx 3 \times 10^{-4},
\end{equation*}
which is much smaller than the prediction accuracy of the neural networks. 
This choice is also motivated by numerical experiments.
\end{rem}

\subsection{Normalisation by a Splitting Approach}
At this stage, the incremental energy density \(\Phi\) depends on \(d + 1\) parameters and shows significant variation over the parameter domain \(\hat{\nu}_{k}\) and \(\alpha_{k}\).  
This pronounced separation between function curves is a challenge for neural networks as it hinders efficient learning.  
Large gaps between function values can prevent smooth interpolation and generalisation, making it difficult to capture underlying patterns during training---unless a large amount of data is used, which becomes intractable even in two spatial dimensions.
Although \(\Phi\) depends only on \(d + 1\) parameters, the parameter \(\hat{\nu}_{k}\) requires a discretisation as fine as the discretisation for the argument \(\hat{\nu}_{k+1}\) since they play a similar role, making the learning computationally infeasible.

To overcome these difficulties, we take advantage of the structure of the pseudo-time incremental energy density function. 
The function \(\Phi\) from \eqref{eq:Phidamage} can be expressed as
\begin{equation*} 
	\Phi(\hat{\nu}_{k+1}; \hat{\nu}_{k}, \alpha_{k})  = \tilde{\Phi}(\hat{\nu}_{k+1}; \alpha_{k}) + \Phi_{\rm shift}(\hat{\nu}_{k}, \alpha_{k}),
\end{equation*}
where \(\tilde{\Phi} \colon \R^{d} \times \R \to \R_{\infty}\) and \(\Phi_{\rm shift}\colon \R^{d} \times \R \to \R_{\infty}\) are defined as
\begin{equation} \label{eq:Phinormalised}
\begin{aligned}
    \tilde{\Phi} (\hat{\nu}_{k+1}; \alpha_{k}) & \coloneqq \varphi(\hat{\nu}_{k+1}, p(\hat{\nu}_{k + 1}; \alpha_{k})) \\
    & \qquad + p(\hat{\nu}_{k + 1}; \alpha_{k})\, D(p(\hat{\nu}_{k + 1}; \alpha_{k})) - \alpha_k \, D(\alpha_k) - \overline{D}(p(\hat{\nu}_{k + 1}; \alpha_{k})) + \overline{D}(\alpha_k)
\end{aligned}
\end{equation}
and 
\begin{equation*}
    \Phi_{\rm shift}(\hat{\nu}_{k}, \alpha_{k}) \coloneqq - \varphi(\hat{\nu}_k,\alpha_k),
\end{equation*}
respectively.
It should be stressed that the function \(\Phi_{\rm shift}\) is independent of \(\hat{\nu}_{k+1}\), hence only dependent on the parameters and constant in the convexification argument \(\hat{\nu}_{k + 1}\). 
Assuming the function \(\varphi^{0}\) is normalised in the sense that \(\inf_{\hat{\nu}} \varphi^{0}(\hat{\nu}) = \varphi^{0}(\mathbbm{1}_d) = 0\), the function \(\tilde{\Phi}\) is also normalised, i.e.~\(\inf_{\hat{\nu}_{k+1}} \tilde{\Phi}(\hat{\nu}_{k+1}; \alpha_{k}) = \tilde{\Phi}(\mathbbm{1}_d) = 0\), where \(\mathbbm{1}_d \in \R^{d}\) denotes the vector containing only ones. 
Consequently, the polyconvex envelope of \(\Phi\) can be obtained from the polyconvexification of the function \(\tilde{\Phi}\) by
\begin{equation} \label{eq:splitpc}
    \Phi^{\pc}(\hat{\nu}_{k+1}; \hat{\nu}_{k}, \alpha_{k}) = \tilde{\Phi}^{\pc}(\hat{\nu}_{k+1}; \alpha_{k}) + \Phi_{\shift}(\hat{\nu}_{k}, \alpha_{k}).
\end{equation}
Therefore, \(\tilde{\Phi}^{\pc}\) should be the focus of an approximation by a neural network or a standard algorithm. Note that this normalisation is domain independent, and just relies on the split \eqref{eq:splitpc}, considering the contribution \(\Phi_{\shift}\) as a shift. 
Removing the dependence on the previous time step \(\hat{\nu}_{k}\) reduces the polyconvexifaction problem to a one-parameter-dependent family, making the learning feasible. 
In what follows, the neural networks are trained to predict the function \(\tilde{\Phi}^{\pc}\) and the function \(\Phi^{\pc}\) is recovered a posteriori by applying the shift \(\Phi_{\rm shift}\) as stated in \eqref{eq:splitpc}.
Since the polyconvex envelopes for both \(\Phi\) and \(\tilde{\Phi}\) are not known analytically, the ground truth is computed according to \cref{sec:SVPC}.

\subsection{Numerical Results for the Two-Dimensional Saint~Venant--Kirchhoff Model}
We consider the function \(\Phi\) from \eqref{eq:Phidamage} in the Saint~Venant--Kirchhoff-based formulation, i.e.~\(\psi^{0}\) from \eqref{eq:STVK} is chosen with the material parameters as before.
We aim for the representation of the function \(\tilde{\Phi}^{\pc}\colon\R^d \times \R_{+} \to \R\), i.e.~the polyconvex envelope of the normalised version \eqref{eq:Phinormalised}, by a neural network. 

\begin{figure}[htbp]
    \centering
    \begin{subfigure}[b]{0.48\textwidth}
    \centering 
    \input{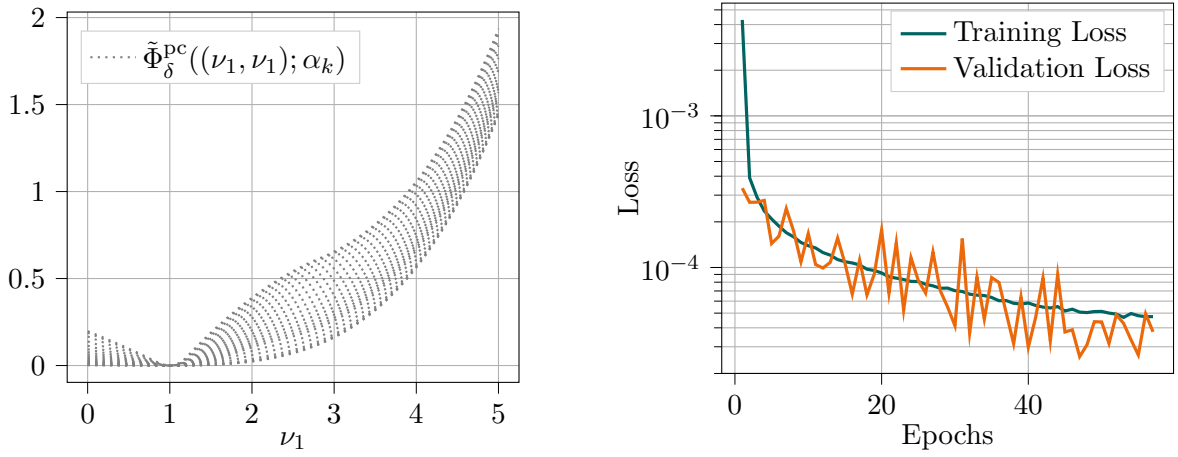}
    \end{subfigure}
    \hfill
    \begin{subfigure}[b]{0.48\textwidth}
    \centering 
    \begin{tikzpicture}
	
	\definecolor{crimson2143940}{RGB}{214,39,40}
	\definecolor{darkgray176}{RGB}{176,176,176}
	\definecolor{darkorange25512714}{RGB}{255,127,14}
	\definecolor{forestgreen4416044}{RGB}{44,160,44}
	\definecolor{lightgray204}{RGB}{204,204,204}
	\definecolor{mediumpurple148103189}{RGB}{148,103,189}
	\definecolor{orchid227119194}{RGB}{227,119,194}
	\definecolor{sienna1408675}{RGB}{140,86,75}
	\definecolor{steelblue31119180}{RGB}{31,119,180}
	
	\begin{axis}[
		width=0.98\textwidth,
		height=6.75cm,
		legend cell align={left},
		legend style={fill opacity=1, draw opacity=1, text opacity=1, draw=lightgray204},
		log basis y={10},
		tick align=outside,
		tick pos=left,
		x grid style={darkgray176},
		xlabel={Epoch},
		xmajorgrids,
		xminorgrids,
		xtick style={color=black},
		y grid style={darkgray176},
		ymajorgrids,
		yminorgrids,
		ymode=log,
		ytick style={color=black},
		xmin=-1.05, xmax=44.05,
		ymin=7.e-07, ymax=0.0852747241544302,
		]
		
		\addplot [very thick, coltrainLoss]
		table {%
			1 0.0473352516097075
			2 0.00019255531399373
			3 9.95738359738893e-05
			4 7.55054310346826e-05
			5 6.47240669133623e-05
			6 5.67101572408843e-05
			7 5.06183915000529e-05
			8 4.74444458743442e-05
			9 4.39515230447996e-05
			10 4.24171420691216e-05
			11 4.04456228877086e-05
			12 3.93035450601192e-05
			13 3.80802008665728e-05
			14 3.72798888229378e-05
			15 3.58670246561466e-05
			16 3.56196458527587e-05
			17 3.44199485384388e-05
			18 3.38954379863712e-05
			19 3.29168822051159e-05
			20 3.19555933639539e-05
			21 3.11163784563531e-05
			22 3.04444169093315e-05
			23 3.03584724318361e-05
			24 2.98393945070364e-05
			25 2.94958337213356e-05
			26 2.92661441353491e-05
			27 2.85225839694521e-05
			28 2.83469110940277e-05
			29 2.85136937189182e-05
			30 2.77761869640375e-05
			31 2.77308711568424e-05
			32 2.76072330586156e-05
			33 2.75785663013519e-05
			34 2.70988295498748e-05
			35 2.68792707897895e-05
			36 2.69184526666668e-05
			37 2.66012289471737e-05
			38 2.63571568312014e-05
			39 2.63454245460509e-05
			40 2.67065136554594e-05
			41 2.58857250000812e-05
			42 2.57295958168123e-05
		};
		\addlegendentry{Train~\(\Loss\)}
				
		\addplot [very thick, colvalLoss]
		table {%
			1 0.000305585566360855
			2 9.79852774869056e-05
			3 7.26124713555525e-05
			4 5.54054141641849e-05
			5 0.000132021247722902
			6 5.31063400122218e-05
			7 6.37770976392225e-05
			8 4.26408588882723e-05
			9 3.59909521114484e-05
			10 4.98531069006887e-05
			11 4.07349788550781e-05
			12 3.18721994288238e-05
			13 3.59983796013185e-05
			14 3.85715604727414e-05
			15 6.57732336940169e-05
			16 3.40316529273425e-05
			17 3.09857846607276e-05
			18 3.52519975350125e-05
			19 2.74517751527814e-05
			20 3.10664148136567e-05
			21 4.19325548131404e-05
			22 4.73680004450309e-05
			23 3.78147720786075e-05
			24 2.50300079342884e-05
			25 2.35684925070429e-05
			26 4.86444251299402e-05
			27 2.65007434815006e-05
			28 2.88506165117194e-05
			29 4.19607635121825e-05
			30 2.52999328928835e-05
			31 2.36927418693276e-05
			32 2.2016563342607e-05
			33 2.65073571096649e-05
			34 2.55989074249278e-05
			35 2.95507913331156e-05
			36 3.83848917434928e-05
			37 2.61682585949552e-05
			38 4.26970142216066e-05
			39 2.76006351535875e-05
			40 2.29259636549566e-05
			41 2.42855628483729e-05
			42 2.53637571593569e-05
		};
		\addlegendentry{Val~\(\Loss\)}
		
		\addplot [very thick, colvalLossMse]
		table {%
			1 0.000169353866528362
			2 6.56559782778426e-05
			3 3.67574708621146e-05
			4 3.34994533804456e-05
			5 3.39364871305998e-05
			6 2.22567511969908e-05
			7 2.80269201779145e-05
			8 2.8104502478758e-05
			9 2.27936597876295e-05
			10 2.45798130074262e-05
			11 1.75247099626323e-05
			12 1.62808388921233e-05
			13 1.5730777916865e-05
			14 2.68909595293341e-05
			15 2.34427081680683e-05
			16 1.68617604571601e-05
			17 2.01124746684611e-05
			18 2.44874909126793e-05
			19 1.31255882189455e-05
			20 1.71668520991914e-05
			21 1.37313677274706e-05
			22 1.44411343307199e-05
			23 1.88734071326784e-05
			24 1.33796681836902e-05
			25 1.39942405357801e-05
			26 1.44199304030689e-05
			27 1.40603537800274e-05
			28 1.19069282555292e-05
			29 1.27393227959409e-05
			30 1.33766487994205e-05
			31 1.10915847470312e-05
			32 1.36025161160414e-05
			33 1.20567029234405e-05
			34 1.40564943922705e-05
			35 2.13100459270217e-05
			36 1.78391077094866e-05
			37 1.91643066677769e-05
			38 2.79326749921212e-05
			39 1.9000851529065e-05
			40 1.22851628482302e-05
			41 1.15123461077204e-05
			42 1.45750404083682e-05
		};
		\addlegendentry{Val~\(\Loss_{\mse}\)}
		
		\addplot [very thick, colvalLossSym]
		table {%
			1 7.12161679879275e-05
			2 1.74476200297771e-05
			3 1.79812016226896e-05
			4 1.3447291546025e-05
			5 4.5810882717962e-05
			6 1.17904915272772e-05
			7 2.08906746944641e-05
			8 1.08837439872732e-05
			9 8.60031443855836e-06
			10 1.80357114169514e-05
			11 1.05258866538645e-05
			12 6.89997318957193e-06
			13 7.3360186390039e-06
			14 7.39478368144142e-06
			15 1.81486383831359e-05
			16 9.67494616111549e-06
			17 7.72245108959167e-06
			18 8.09288103963312e-06
			19 6.24083665839978e-06
			20 9.17175813649275e-06
			21 9.83958152139799e-06
			22 1.25716694294885e-05
			23 1.37500652372688e-05
			24 6.40369669549187e-06
			25 5.61715902185633e-06
			26 1.72853566727074e-05
			27 7.38176633168214e-06
			28 6.9869178224051e-06
			29 1.11305834595656e-05
			30 6.98043348971098e-06
			31 5.08015371752922e-06
			32 4.92744539835468e-06
			33 7.40215442078597e-06
			34 7.55945571410369e-06
			35 6.32356939404323e-06
			36 1.34356653246872e-05
			37 5.31103645469245e-06
			38 1.33037368567672e-05
			39 6.28839429138859e-06
			40 6.2927197761773e-06
			41 7.15614293315342e-06
			42 6.61794471538975e-06
		};
		\addlegendentry{Val~\(\lambda_{\sym} \, \Loss_{\sym}\)}
	
		\addplot [very thick, colvalLossIneq]
		table {%
			1 6.50155319429613e-05
			2 1.48816792529998e-05
			3 1.78737988339062e-05
			4 8.45866921805851e-06
			5 5.22738778462023e-05
			6 1.9059097312083e-05
			7 1.48595027499636e-05
			8 3.65261246777956e-06
			9 4.5969778913109e-06
			10 7.23758249366205e-06
			11 1.26843822483847e-05
			12 8.69138733094746e-06
			13 1.29315830945921e-05
			14 4.28581727936316e-06
			15 2.41818871013166e-05
			16 7.49494633316753e-06
			17 3.15085891655307e-06
			18 2.67162555646104e-06
			19 8.08535026255125e-06
			20 4.72780458289176e-06
			21 1.83616055809671e-05
			22 2.03551966715954e-05
			23 5.19129972698253e-06
			24 5.2466430771805e-06
			25 3.9570929588167e-06
			26 1.69391380223161e-05
			27 5.05862337381138e-06
			28 9.95677043721716e-06
			29 1.80908572742555e-05
			30 4.94285061928761e-06
			31 7.52100341413493e-06
			32 3.48660182530369e-06
			33 7.04849976055151e-06
			34 3.98295731638803e-06
			35 1.91717599812891e-06
			36 7.11011868028905e-06
			37 1.69291545823286e-06
			38 1.46060241885673e-06
			39 2.31138934559947e-06
			40 4.34808102590082e-06
			41 5.61707382810372e-06
			42 4.17077204073973e-06
		};
		\addlegendentry{Val~\(\lambda_{\ineq} \, \Loss_{\ineq}\)}

	\end{axis}
\end{tikzpicture}
    \end{subfigure}
    \caption{
    Left: Learning data for STVK based damage.
    Diagonal cross-sections of the target values, i.e.~polyconvex envelopes \(\tilde{\Phi}^{\pc}_{\delta}\) for the Saint~Venant--Kirchhoff-based function \(\tilde{\Phi}\) with $\alpha_k \in \mathcal{I}_{\alpha_k}$.
    Right: Learning curves for STVK based damage.
    Learning curves based on \(\Loss\) from \eqref{eq:loss} for a single network initialisation for the Saint~Venant--Kirchhoff-based function \(\tilde{\Phi}^{\pc}_{\delta}\).
    }
    \label{fig:Training_Loss_damage_STVK}
\end{figure}

\subsubsection{Network Architecture and Learning Data}
In this example, we implement a PICNN consisting of three hidden layers, where the \(\hat{m}\)-path consists of layers with \(30\), \(60\) and \(60\) neurons and the \(\zeta\)-path, i.e.~the parameter path, consists of \(15\), \(30\) and \(30\) neurons, respectively. 
The convex input \(\hat{m}\) represents the minors of the signed singular values \(\hat{\nu}_{k+1}\), i.e.~\(\m(\hat{\nu}_{k+1})\), while the nonconvex inputs, denoted as \(\zeta\), correspond to the parameter \(\alpha_k\). 
The learning domain for the parameter \(\alpha_k\) is set to \([0, \alpha_\infty]\) while the learning domain for \(\hat{\nu}_{k+1}\) is defined as \([\nu_{\min},\nu_{\max}]^2\) with \(\nu_{\min}=0.1\) and \(\nu_{\max}=5\) and all permutations included in \(\Pid\).
The data points in the singed singular value space are obtained by discretising each axis on the interval \([\nu_{\min}, \nu_{\max}]\) by \([\nu_{\min} : 0.005 : 1.2] \cup [1.2 : 0.02 : \nu_{\max}]\), leading to a \(411 \times 411\) grid in \(\R_{+}^d\). 
Additionally all possible permutations of these points induced by transformations in \(\Pid\) are included to extend the data to full \(\R^d\).
Since \(\tilde{\Phi}^{\pc}(\hat{\nu}; \alpha_{k}) = + \infty\), for \(\nu_1 \, \nu_2 \leq 0\), only the quadrants with positive product \(\nu_1 \, \nu_2\) are considered in the learning data set, resulting in \num{337842} grid points in the signed singular value space. 
Notably, the learning domain for \(\hat{\nu}_{k+1}\) covers the part of the signed singular value space associated to positive determinant.
For \(\alpha_k\), we choose the set of learning values as \(\mathcal{I}_{\alpha_k} = \{0, 0.1, \dotsc, 1.5, 1.75, 2, 2.25, 2.5, 3, 3.5, 4\}\), i.e.~23 values. 
These discretisations lead to a learning dataset consisting of \num{7770366} points. 
The target values \(\tilde{\Phi}^{\pc}_{\delta}\) are computed with the SVPC~LP algorithm described in \cref{sec:SVPC}, with discretisation radius \(r = 5.1\).
The lattice width is varied, since points close to the origin require a finer resolution of the computational grid to resolve the growth of the function towards the determinant constrained regime. 
That is why for \(\hat{\nu}\) with \(\min(\nu_1, \nu_2) \leq 0.125\), a lattice size of \(\delta \approx 0.01\), for \(\hat{\nu}\) with \(\min(\nu_1, \nu_2) \leq 1.25\), a lattice size \(\delta \approx 0.02\) and for all other evaluation points \(\delta = 0.04\) is chosen. For this discretisation strategy the computation of the target values takes approximately \num{400} CPU hours.  \cref{fig:Training_Loss_damage_STVK} illustrates the target values of the learning data along the positive part of the diagonal axis \((\nu_1, \nu_1)\). 
From the learning dataset, \(70 \%\) of the samples are randomly assigned to the training set, while the remaining \(30 \%\) constitute the validation set. 
For this experiment the patience parameter is set to \(10\) and the penalty parameters included in the loss function are set to \(\lambda_{\ineq}=4\) and \(\lambda_{\sym}=2\).

\subsubsection{Numerical Results}
On average, for a single realisation, the training process requires \(1.5\)~h \(\pm\) \(18\) minutes to complete \(38 \pm 8\) epochs. 
The final training loss is \(3.11 \times 10^{-5}\pm 4.9 \times 10^{-6}\), while the validation loss reaches \(3.22  \times 10^{-5} \pm 6.6 \times 10^{-6}\). 
The learning curves for one of these realisations are presented in \cref{fig:Training_Loss_damage_STVK}, including the individual contributions to the validation loss evaluation. 

\begin{figure}[htbp]
	\centering
	\begin{subfigure}{0.49\textwidth}
		\centering
		\begin{tikzpicture}
	\begin{axis}[
		width=0.84\textwidth,
		height=0.58\textwidth,
		scale only axis,
		legend pos=north east,
		xmin=-0.25, xmax=7.25,
		ymin=0.0, ymax=0.01, 
		scaled y ticks=false,
		yticklabel style={
			/pgf/number format/fixed,
			/pgf/number format/precision=3
		},
		grid=both,	
		major grid style={line width=0.4pt, draw=gray!60},
		minor grid style={line width=0.2pt, draw=gray!30},
		xlabel={\(\alpha_{k}\)},
		minor tick num=1,
		xtick={0.,0.5,1.,1.5, 2, 2.5, 3, 3.5, 4, 5, 6, 7},
		xticklabels={0,,1,, 2,, 3,, 4, 5, 6, 7},
		]

		\addplot+[
			unia-green,
			only marks,
			mark=*,
			mark size=2.5pt,
			draw=unia-green,
			mark options={fill=unia-green},
			error bars/.cd,
				y dir=both,
				y explicit,
				error bar style={line width=0.8pt, color=unia-green},
		] table[x=alpha, y=mean, y error=std] {Pictures/RealDamage2D/STVK/NN4poly_Iso_damage_STVK_mean_errs.dat};
		
		\addlegendentry{Avg mean errs \(\pm\) std dev}
		
	\end{axis}
\end{tikzpicture}
		\caption{Average of mean errors}
	\end{subfigure}
	\hfill
	\begin{subfigure}{0.49\textwidth}
		\centering
		\begin{tikzpicture}
	\begin{axis}[
		width=0.84\textwidth,
		height=0.58\textwidth,
		scale only axis,
		legend pos=north east,
		xmin=-0.25, xmax=7.25,
		ymin=0.0, ymax=0.025, 
		scaled y ticks=false,
		yticklabel style={
			/pgf/number format/fixed,
			/pgf/number format/precision=3
		},
		grid=both,	
		major grid style={line width=0.4pt, draw=gray!60},
		minor grid style={line width=0.2pt, draw=gray!30},
		xlabel={\(\alpha_{k}\)},
		minor tick num=1,
		xtick={0.,0.5,1.,1.5, 2, 2.5, 3, 3.5, 4, 5, 6, 7},
		xticklabels={0,,1,, 2,, 3,, 4, 5, 6, 7},
		]
		
		\addplot+[
			unia-green,
			only marks,
			mark=*,
			mark size=2.5pt,
			draw=unia-green,
			mark options={fill=unia-green},
			error bars/.cd,
			y dir=both,
			y explicit,
			error bar style={line width=0.8pt, color=unia-green},
		] table[x=alpha, y=mean, y error=std] {Pictures/RealDamage2D/STVK/NN4poly_Iso_damage_STVK_rel_quad_errs.dat};
		
		\addlegendentry{Avg rel quad errs \(\pm\) std dev}
		
	\end{axis}
\end{tikzpicture}
		\caption{Average of relative quadratic errors}
	\end{subfigure}	
	\begin{subfigure}{0.49\textwidth}
		\centering
		\begin{tikzpicture}
	\begin{axis}[
		width=0.84\textwidth,
		height=0.58\textwidth, 
		scale only axis,
		legend pos=north east,
		xmin=-0.25, xmax=7.25,
		ymin=0.0, ymax=0.03, 
		scaled y ticks=false,
		yticklabel style={
			/pgf/number format/fixed,
			/pgf/number format/precision=3
		},
		grid=both,	
		major grid style={line width=0.4pt, draw=gray!60},
		minor grid style={line width=0.2pt, draw=gray!30},
		xlabel={\(\alpha_{k}\)},
		minor tick num=1,
		xtick={0.,0.5,1.,1.5, 2, 2.5, 3, 3.5, 4, 5, 6, 7},
		xticklabels={0,,1,, 2,, 3,, 4, 5, 6, 7},
		]
		
		\addplot+[
			unia-green,
			only marks,
			mark=*,
			mark size=2.5pt,
			draw=unia-green,
			mark options={fill=unia-green},
			error bars/.cd,
			y dir=both,
			y explicit,
			error bar style={line width=0.8pt, color=unia-green},	
		] table[x=alpha, y=mean, y error=std] {Pictures/RealDamage2D/STVK/NN4poly_Iso_damage_STVK_rel_max_errs.dat};
		
		\addlegendentry{Avg rel max errs \(\pm\) std dev}
		
	\end{axis}
\end{tikzpicture}
		\caption{Average of relative max errors}
	\end{subfigure}
	\caption{
		Prediction errors of \(\tilde{\Phi}^\pc_\pred\) (average over ten network realisations) with respect to the reference polyconvex envelopes \(\tilde{\Phi}^\pc_\delta\) computed via SVPC~LP algorithm for different values \(\alpha_k\) for the Saint~Venant--Kirchhoff-based damage model function \(\tilde{\Phi}\).
	}
	\label{fig:Errors_STVK}
\end{figure}

\Cref{fig:Errors_STVK} illustrates the relative errors between the predictions and the polyconvex envelopes computed by the algorithm from \cref{sec:SVPC} for different values of \(\alpha_k\), including values that are outside the training domain.  
These errors are computed on a uniform \(100 \times 100\) lattice of the domain \([0.1, 5]^2 \subset [\nu_{\min},\nu_{\max}]^2\). 
Across all considered values \(\alpha_k\), the errors remain consistently low, ranging between 1\% and 2\%. 
In particular, it is important to note that the hypothesis stating that the energies become independent of \(\alpha_k\) for \(\alpha_k\geq \alpha_{\infty}\) is verified and validates the choice of \(\alpha_\infty=4\).

\begin{figure}[htbp]
	\centering	
	\input{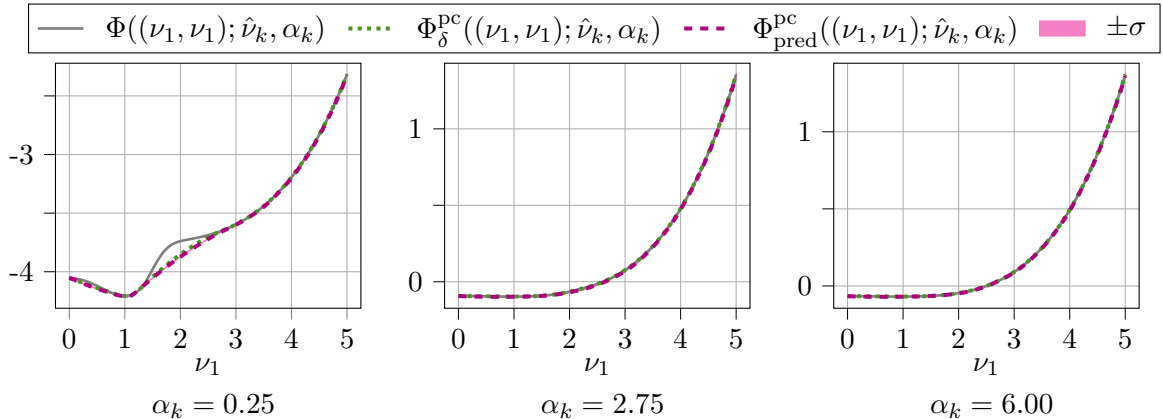}
	\caption{
		Comparison of the polyconvex envelopes \(\Phi^{\pc}_\pred\) (averaged over ten network realisations) and \(\Phi^{\pc}_{\delta}\) for the Saint~Venant--Kirchhoff based damage model with \(\hat{\nu}_k=[2.5,2.5]\) and for \(\alpha_k \in \{0.25, 1.25\}\) on \num{200} points on \(\nu_1 \in [0.1, 5.0]\). 
		The plot shows the cross-section along the diagonal axis. The standard deviation of the ten predictions is marked by \(\sigma\).
	}
	\label{fig:STVK_curves}
\end{figure}

Additionally, \cref{fig:STVK_curves} provides examples of re-shifted predictions, following \eqref{eq:splitpc}, which are in good agreement with the polyconvex envelopes computed by the SVPC~LP algorithm. 
In particular, these predictions are obtained for values of \(\alpha_k\) that are not included in the learning set \(\mathcal{I}_{\alpha_k}\) and lie outside the parameter's training domain, demonstrating good predictability and generalisation capabilities of the network.

For comparison, the computation of one polyconvex envelope on a \(100 \times 100\)-grid on the box \([0.1, 5]^2\) using the SVPC~LP algorithm of \cref{sec:SVPC}, with the same discretisation parameters as above, takes \(49\) minutes on a single CPU (exploiting the symmetry due to \(\Pid\)-invariance) while its prediction via the neural network takes only \(0.05\) seconds, emphasising once again the benefits of using a neural-network compression approach for the polyconvexification in engineering applications.
In addition, the neural network implemented in this example contains \num{14800} parameters, comprising both weights and biases. 
Once again, this number of parameters remains smaller than the storage required for only two polyconvex envelopes on a \(100 \times 100\) grid.

\subsection{Numerical Results for the Two-Dimensional Neo-Hookean Model}
We consider the function \(\Phi\) from \eqref{eq:Phidamage} in the neo-Hookean based formulation, i.e.~\(\psi^{0}\) from \eqref{eq:NH} is chosen and the material parameters are set as before.
We aim for the representation of the function \(\tilde{\Phi}^{\pc}\colon\R^d \times \R_{+} \to \R\), i.e.~the polyconvex envelope of the normalised version \eqref{eq:Phinormalised}, by a neural network. 

\begin{figure}[htbp]
    \centering
    \begin{subfigure}[b]{0.48\textwidth}
    \centering 
    \input{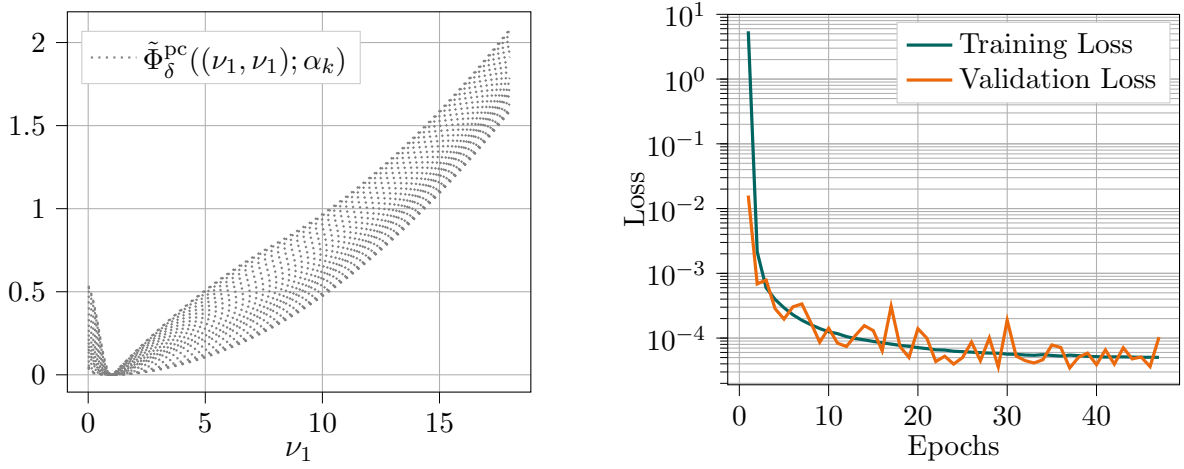}
    \end{subfigure}
    \hfill
    \begin{subfigure}[b]{0.48\textwidth}
    \centering 
\begin{tikzpicture}
	
	\definecolor{darkgray176}{RGB}{176,176,176}
	\definecolor{darkorange}{RGB}{255,140,0}
	\definecolor{lightgray204}{RGB}{204,204,204}
	\definecolor{royalblue}{RGB}{65,105,225}
	
	\begin{axis}[
		width=0.98\textwidth,
		legend cell align={left},
		legend style={fill opacity=1, draw opacity=1, text opacity=1, draw=lightgray204},
		log basis y={10},
		tick align=outside,
		tick pos=left,
		x grid style={darkgray176},
		xlabel={Epoch},
		xmajorgrids,
		xminorgrids,
		xtick style={color=black},
		y grid style={darkgray176},
		ymajorgrids,
		yminorgrids,
		ymode=log,
		ytick style={color=black},		
		xmin=-2.15, xmax=67.15,
		ymin=1.06653207844255e-07, ymax=104.672514829983,
		]
		
		\addplot [very thick, coltrainLoss]
		table {%
			1 40.8424519588944
			2 0.00226933934938455
			3 0.000313898310884045
			4 0.000180864841738203
			5 8.90304086561113e-05
			6 6.82094532197758e-05
			7 6.15042180940062e-05
			8 6.04445597359317e-05
			9 5.33451996418096e-05
			10 5.00230007239492e-05
			11 4.56863770966805e-05
			12 4.68364326728978e-05
			13 4.25050582874227e-05
			14 4.07808969549715e-05
			15 3.94882675252529e-05
			16 3.80312156567729e-05
			17 3.73150995579804e-05
			18 3.56148683796647e-05
			19 3.467655919987e-05
			20 3.48931581080423e-05
			21 3.3857425022585e-05
			22 3.34724116742393e-05
			23 3.28479461595217e-05
			24 3.27195271282632e-05
			25 3.22136080230992e-05
			26 3.20692240634449e-05
			27 3.16018682584157e-05
			28 3.13643137337259e-05
			29 3.12081885249704e-05
			30 3.09009232429736e-05
			31 3.07034220273585e-05
			32 3.03500571917796e-05
			33 3.02742578983933e-05
			34 3.21221410315209e-05
			35 2.95656150249401e-05
			36 2.98096834229037e-05
			37 2.96634141961044e-05
			38 2.91591632342259e-05
			39 2.91667100859575e-05
			40 2.9083628855556e-05
			41 2.88198713601338e-05
			42 2.88430655093407e-05
			43 2.88221516690396e-05
			44 2.84939932272295e-05
			45 2.83373531324109e-05
			46 2.93902113350711e-05
			47 2.79347496828271e-05
			48 2.83628149548513e-05
			49 2.75272479222536e-05
			50 2.73935369296628e-05
			51 2.71766204894408e-05
			52 2.7835578260826e-05
			53 2.75214565374091e-05
			54 2.6840329234386e-05
			55 2.79759727729765e-05
			56 2.65492705583524e-05
			57 2.63918421086953e-05
			58 2.65949336090428e-05
			59 2.61267429014265e-05
			60 2.6638474645593e-05
			61 2.57625743902301e-05
			62 2.57843401651251e-05
			63 2.57219531347035e-05
			64 2.62007746634186e-05
		};
		\addlegendentry{Train~\(\Loss\)}

		\addplot [very thick, colvalLoss]
		table {%
			1 0.000863218228596513
			2 0.000978234317700199
			3 0.000119484735904885
			4 9.72064491898933e-05
			5 6.71061189306294e-05
			6 5.63456400264961e-05
			7 5.45449178769474e-05
			8 4.79586176547765e-05
			9 4.24673506167265e-05
			10 5.77797396496256e-05
			11 7.52654881898462e-05
			12 5.82991624305799e-05
			13 4.20425351548088e-05
			14 3.23834404312674e-05
			15 4.23036088811906e-05
			16 3.37101846121949e-05
			17 3.32417460082268e-05
			18 3.479835746149e-05
			19 3.04210642861311e-05
			20 3.05482802615947e-05
			21 3.01283757924714e-05
			22 3.14060676523287e-05
			23 2.90082350906785e-05
			24 3.08402889024319e-05
			25 3.434327528719e-05
			26 4.06821402079213e-05
			27 3.52776760261853e-05
			28 2.79988469746904e-05
			29 4.82107191227402e-05
			30 2.71827547630272e-05
			31 2.9892904413261e-05
			32 2.78259117798729e-05
			33 4.91126410182111e-05
			34 2.57749305205607e-05
			35 2.88253526138536e-05
			36 3.59980294815106e-05
			37 2.75801594605022e-05
			38 4.88124667507322e-05
			39 2.39636298484408e-05
			40 2.42121602270375e-05
			41 2.47850049044634e-05
			42 2.7145575416986e-05
			43 2.73806348877481e-05
			44 3.19958563184086e-05
			45 2.54344832380898e-05
			46 2.37700916880167e-05
			47 3.07398056181857e-05
			48 2.62190339588581e-05
			49 2.27564877950615e-05
			50 2.30050258626949e-05
			51 2.42750989326285e-05
			52 2.6267010628808e-05
			53 2.27269641989975e-05
			54 2.19671853415549e-05
			55 2.40325145387688e-05
			56 2.26422413795223e-05
			57 2.38536409238694e-05
			58 2.39464403336221e-05
			59 2.37354389591035e-05
			60 2.23155065630854e-05
			61 4.25274939940524e-05
			62 3.35813669747443e-05
			63 2.25656814098597e-05
			64 2.83978870674762e-05
		};
		\addlegendentry{Val~\(\Loss\)}

		\addplot [very thick, colvalLossMse]
		table {%
			1 0.000504362175034665
			2 0.000581963272909111
			3 6.83125702497883e-05
			4 4.58225431093843e-05
			5 3.63750326275494e-05
			6 3.81338201644515e-05
			7 3.47503314782552e-05
			8 2.32921041083552e-05
			9 2.09473646966461e-05
			10 3.60819883163164e-05
			11 2.84420575638961e-05
			12 1.97702639541198e-05
			13 1.67800791105475e-05
			14 1.95382697415144e-05
			15 2.94681475396506e-05
			16 1.56236386566906e-05
			17 1.90552372957884e-05
			18 2.60320921961496e-05
			19 1.59533305810873e-05
			20 1.66520867393363e-05
			21 1.83096266568566e-05
			22 2.17812901463632e-05
			23 1.8586791337353e-05
			24 1.29841899274198e-05
			25 2.39485694113887e-05
			26 2.55697974215919e-05
			27 2.17853897393398e-05
			28 1.32768290382978e-05
			29 1.76544147202975e-05
			30 1.67668651225506e-05
			31 1.45619952601486e-05
			32 1.72922099269754e-05
			33 1.42195511492938e-05
			34 1.66113548475594e-05
			35 1.42007691443766e-05
			36 1.29065299878478e-05
			37 1.69837912296555e-05
			38 2.75666360933343e-05
			39 1.45548839011533e-05
			40 1.45748306007521e-05
			41 1.48438823043821e-05
			42 1.42696348232162e-05
			43 1.38496062605953e-05
			44 2.21388750521311e-05
			45 1.52642444631251e-05
			46 1.24961693652564e-05
			47 1.7820566461603e-05
			48 1.32292614416824e-05
			49 1.26181241738226e-05
			50 1.30437199176879e-05
			51 1.31904478670371e-05
			52 1.55113680259791e-05
			53 1.3112263386808e-05
			54 1.16666849502345e-05
			55 1.20285378561792e-05
			56 1.51498689708722e-05
			57 1.34521831098582e-05
			58 1.32469133587007e-05
			59 1.59493721742009e-05
			60 1.2766666924789e-05
			61 1.31815358221391e-05
			62 2.62558342614847e-05
			63 1.29855211616255e-05
			64 2.13193219946114e-05
		};
		\addlegendentry{Val~\(\Loss_{\mse}\)}

		\addplot [very thick, colvalLossSym]
		table {%
			1 0.000172542467712624
			2 0.000252956409571282
			3 2.65223842528402e-05
			4 3.01564773889434e-05
			5 1.7292676337627e-05
			6 1.11592033390504e-05
			7 1.24392460088346e-05
			8 1.14188430755823e-05
			9 9.78235483112541e-06
			10 1.80097864051961e-05
			11 1.2002987656663e-05
			12 1.01317960262637e-05
			13 9.84005247878279e-06
			14 8.10706509581981e-06
			15 9.87287742833556e-06
			16 6.92043283014736e-06
			17 9.695167605196e-06
			18 6.87082115770147e-06
			19 7.69032399368281e-06
			20 7.33710447786943e-06
			21 7.82907131462875e-06
			22 7.04427538779746e-06
			23 7.10772829147682e-06
			24 6.57908744568066e-06
			25 8.69212800666812e-06
			26 1.25645541041581e-05
			27 1.09084497378227e-05
			28 7.18829606748997e-06
			29 8.00897633213581e-06
			30 6.71097540691533e-06
			31 8.83657507890144e-06
			32 7.83117928772771e-06
			33 1.08463342851768e-05
			34 6.34439445774206e-06
			35 8.76913916912114e-06
			36 1.0687225654373e-05
			37 7.46435866186041e-06
			38 1.8996619916612e-05
			39 6.3313336680599e-06
			40 6.40971318386234e-06
			41 6.3876613356377e-06
			42 7.29170133719373e-06
			43 8.65427676360599e-06
			44 8.30676794309824e-06
			45 6.5873253612809e-06
			46 6.58964797887023e-06
			47 9.16299455867542e-06
			48 8.2502732897202e-06
			49 6.09095474665198e-06
			50 6.4543219483652e-06
			51 6.41666516168753e-06
			52 7.95940629844536e-06
			53 6.34303619317782e-06
			54 5.8453043167773e-06
			55 6.63427571475352e-06
			56 5.69705416127869e-06
			57 7.12304316005462e-06
			58 6.98284638280685e-06
			59 5.94582513030412e-06
			60 6.27852439050303e-06
			61 1.33196497125078e-05
			62 6.15969825390898e-06
			63 6.079556702647e-06
			64 5.98522627189723e-06
		};
		\addlegendentry{Val~\(\lambda_{\sym} \, \Loss_{\sym}\)}

		\addplot [very thick, colvalLossIneq]
		table {%
			1 0.00018631358630997
			2 0.000143314635691023
			3 2.46497812923631e-05
			4 2.12274287631037e-05
			5 1.34384099504675e-05
			6 7.05261649156536e-06
			7 7.35534037139386e-06
			8 1.324767050906e-05
			9 1.17376310675247e-05
			10 3.68796497438316e-06
			11 3.48204430055253e-05
			12 2.83971024240513e-05
			13 1.54224035527131e-05
			14 4.73810560296163e-06
			15 2.96258395006232e-06
			16 1.11661131353528e-05
			17 4.49134112316159e-06
			18 1.89544409585351e-06
			19 6.77740972569192e-06
			20 6.55908906537475e-06
			21 3.98967784541363e-06
			22 2.58050212428071e-06
			23 3.31371546049472e-06
			24 1.12770115566001e-05
			25 1.70257787528856e-06
			26 2.54778867498302e-06
			27 2.58383656797783e-06
			28 7.53372187156466e-06
			29 2.25473280819137e-05
			30 3.70491423639432e-06
			31 6.49433408719731e-06
			32 2.70252256830808e-06
			33 2.40467555691971e-05
			34 2.81918122649836e-06
			35 5.85544430032089e-06
			36 1.24042738435198e-05
			37 3.1320095945306e-06
			38 2.2492107212537e-06
			39 3.07741227657252e-06
			40 3.22761644638224e-06
			41 3.55346127505511e-06
			42 5.58423927140531e-06
			43 4.87675186197871e-06
			44 1.55021331476135e-06
			45 3.58291340507689e-06
			46 4.68427434887438e-06
			47 3.75624459778993e-06
			48 4.73949922547782e-06
			49 4.04740888699267e-06
			50 3.50698399389231e-06
			51 4.66798589363321e-06
			52 2.79623627597572e-06
			53 3.27166460754359e-06
			54 4.45519607314218e-06
			55 5.36970094094442e-06
			56 1.79531826826949e-06
			57 3.27841463190127e-06
			58 3.71668059781304e-06
			59 1.84024164042864e-06
			60 3.2703152371525e-06
			61 1.60263084486877e-05
			62 1.16583445824379e-06
			63 3.50060354670275e-06
			64 1.09333881236892e-06
		};
		\addlegendentry{Val~\(\lambda_{\ineq} \, \Loss_{\ineq}\)}
	\end{axis}
	
\end{tikzpicture}
    \end{subfigure}
    \caption{
    Left:
    Learning data for NH based damage.
    Diagonal cross-sections of the target values, i.e.~polyconvex envelopes \(\tilde{\Phi}^{\pc}_{\delta}\) for the neo-Hookean-based function \(\tilde{\Phi}\) with $\alpha_k \in \mathcal{I}_{\alpha_k}$.
    Right: Learning curves for NH-based damage.
    Learning curves based on \(\Loss\) from \eqref{eq:loss} for a single network initialisation for the neo-Hookean-based function \(\tilde{\Phi}^{\pc}_{\delta}\).
    } 
    \label{fig:Training_Loss_damage_NH}
\end{figure}

\subsubsection{Network Architecture and Learning Data}
As before, we implement a PICNN consisting of three hidden layers, where the \(\hat{m}\)-path consists of layers with \(30\), \(60\) and \(60\) neurons and the \(\zeta\)-path, i.e.~the parameter path, consists of \(15\), \(30\) and \(30\) neurons, respectively. 
The convex input \(\hat{m}\) represents the minors of the signed singular values \(\hat{\nu}_{k+1}\), i.e.~\(\m(\hat{\nu}_{k+1})\), while the nonconvex inputs, denoted as \(\zeta\), correspond to the parameter \(\alpha_k\). 
The learning domain for the parameter \(\alpha_k\) is set to \([0, \alpha_\infty]\) while the learning domain for \(\hat{\nu}_{k+1}\) is defined as \([\nu_{\min},\nu_{\max}]^2\) with \(\nu_{\min}=0.55\) and \(\nu_{\max}=18\) and all permutations included in \(\Pid\).
The data points in the singed singular value space are obtained by discretising each axis on the interval \([\nu_{\min}, \nu_{\max}]\) by \([\nu_{\min} : 0.015 : 1.75] \cup [1.75 : 0.05 : 18]\), leading to a \(406 \times 406\) grid in \(\R_{+}^d\). 
Additionally, all possible permutations of these points induced by transformations in the symmetry group \(\Pid\) are included to extend the data to full \(\R^d\).
Since \(\tilde{\Phi}^{\pc}(\hat{\nu}; \alpha_{k}) = + \infty\), for \(\nu_1 \, \nu_2 \leq 0\), only the quadrants with positive product \(\nu_1 \, \nu_2\) are considered in the learning data set, resulting in \num{329672} grid points in the signed singular value space. 
For \(\alpha_k\), we choose the set of learning values as \(\mathcal{I}_{\alpha_k} = \{0, 0.1, \dotsc, 1.5, 1.75, 2, 2.25, 2.5, 3, 3.5, 4\}\), i.e.~23 values. 
These discretisations lead to a learning dataset consisting of \num{7582456} points. 

The target values \(\tilde{\Phi}^{\pc}_{\delta}\) are computed with the SVPC~LP algorithm described in \cref{sec:SVPC}, with discretisation radius \(r = 20\).
The lattice width is varied, since points close to the origin require a finer resolution of the computational grid to resolve the growth of the function towards the determinant-constrained regime. 
That is why for \(\hat{\nu}\) with \(\min(\nu_1, \nu_2) \leq 1.25\), a lattice size \(\delta \approx 0.039\) and for all other evaluation points \(\delta = 0.078\) is chosen. For this discretisation strategy the computation of the target values takes approximately \num{600} CPU hours. 
\cref{fig:Training_Loss_damage_NH} illustrates the target values of the learning data along the positive part of the diagonal axis \((\nu_1, \nu_1)\). 
From the learning dataset, \(70 \%\) of the samples are randomly assigned to the training set, while the remaining \(30 \%\) constitute the validation set. 
For this experiment, the patience parameter is set to \(10\) and the penalty parameters included in the loss function are set to \(\lambda_{\ineq}=4\) and \(\lambda_{\sym}=2\).

\subsubsection{Numerical Results}
On average, for a single realisation, the training process requires \(3.2 \pm 1.2\)~h to complete \(63 \pm 22\) epochs. 
The final training loss is \(3.04 \times 10^{-5}\pm 8.2 \times 10^{-6}\), while the validation loss reaches \(3.57 \times 10^{-5} \pm 1.6 \times 10^{-5}\). 
The learning curves for one of these realisations are presented in \cref{fig:Training_Loss_damage_NH}.

\begin{figure}[htbp]
	\centering
	\begin{subfigure}{0.49\textwidth}
		\centering
		\begin{tikzpicture}
	\begin{axis}[
		width=0.84\textwidth,
		height=0.58\textwidth,
		scale only axis,
		legend pos=north east,
		xmin=-0.25, xmax=7.25,
		ymin=0.0, ymax=0.01, 
		scaled y ticks=false,
		yticklabel style={
			/pgf/number format/fixed,
			/pgf/number format/precision=3
		},
		grid=both,	
		major grid style={line width=0.4pt, draw=gray!60},
		minor grid style={line width=0.2pt, draw=gray!30},
		xlabel={\(\alpha_{k}\)},
		minor tick num=1,
		xtick={0.,0.5,1.,1.5, 2, 2.5, 3, 3.5, 4, 5, 6, 7},
		xticklabels={0,,1,, 2,, 3,, 4, 5, 6, 7},
		]

		\addplot+[
			unia-green,
			only marks,
			mark=*,
			mark size=2.5pt,
			draw=unia-green,
			mark options={fill=unia-green},
			error bars/.cd,
				y dir=both,
				y explicit,
				error bar style={line width=0.8pt, color=unia-green},
		] table[x=alpha, y=mean, y error=std] {Pictures/RealDamage2D/NH/NN4poly_Iso_damage_NH_mean_errs.dat};
		
		\addlegendentry{Avg mean errs \(\pm\) std dev}
		
	\end{axis}
\end{tikzpicture}
		\caption{Average of mean errors}
	\end{subfigure}
	\hfill
	\begin{subfigure}{0.49\textwidth}
		\centering
		\begin{tikzpicture}
	\begin{axis}[
		width=0.84\textwidth,
		height=0.58\textwidth,
		scale only axis,
		legend pos=north east,
		xmin=-0.25, xmax=7.25,
		ymin=0.0, ymax=0.025, 
		scaled y ticks=false,
		yticklabel style={
			/pgf/number format/fixed,
			/pgf/number format/precision=3
		},
		grid=both,	
		major grid style={line width=0.4pt, draw=gray!60},
		minor grid style={line width=0.2pt, draw=gray!30},
		xlabel={\(\alpha_{k}\)},
		minor tick num=1,
		xtick={0.,0.5,1.,1.5, 2, 2.5, 3, 3.5, 4, 5, 6, 7},
		xticklabels={0,,1,, 2,, 3,, 4, 5, 6, 7},
		]
		
		\addplot+[
			unia-green,
			only marks,
			mark=*,
			mark size=2.5pt,
			draw=unia-green,
			mark options={fill=unia-green},
			error bars/.cd,
			y dir=both,
			y explicit,
			error bar style={line width=0.8pt, color=unia-green},
		] table[x=alpha, y=mean, y error=std] {Pictures/RealDamage2D/NH/NN4poly_Iso_damage_NH_rel_quad_errs.dat};
		
		\addlegendentry{Avg rel quad errs \(\pm\) std dev}
		
	\end{axis}
\end{tikzpicture}
		\caption{Average of relative quadratic errors}
	\end{subfigure}	
	\begin{subfigure}{0.49\textwidth}
		\centering
		\begin{tikzpicture}
	\begin{axis}[
		width=0.84\textwidth,
		height=0.58\textwidth, 
		scale only axis,
		legend pos=north east,
		xmin=-0.25, xmax=7.25,
		ymin=0.0, ymax=0.035, 
		scaled y ticks=false,
		yticklabel style={
			/pgf/number format/fixed,
			/pgf/number format/precision=3
		},
		grid=both,	
		major grid style={line width=0.4pt, draw=gray!60},
		minor grid style={line width=0.2pt, draw=gray!30},
		xlabel={\(\alpha_{k}\)},
		minor tick num=1,
		xtick={0.,0.5,1.,1.5, 2, 2.5, 3, 3.5, 4, 5, 6, 7},
		xticklabels={0,,1,, 2,, 3,, 4, 5, 6, 7},
		]
		
		\addplot+[
			unia-green,
			only marks,
			mark=*,
			mark size=2.5pt,
			draw=unia-green,
			mark options={fill=unia-green},
			error bars/.cd,
			y dir=both,
			y explicit,
			error bar style={line width=0.8pt, color=unia-green},	
		] table[x=alpha, y=mean, y error=std] {Pictures/RealDamage2D/NH/NN4poly_Iso_damage_NH_rel_max_errs.dat};
		
		\addlegendentry{Avg rel max errs \(\pm\) std dev}
		
	\end{axis}
\end{tikzpicture}
		\caption{Average of relative max errors}
	\end{subfigure}
    \caption{
		Prediction errors of \(\tilde{\Phi}^\pc_\pred\) (average over ten network realisations) with respect to the reference polyconvex envelopes \(\tilde{\Phi}^\pc_\delta\) computed via SVPC~LP algorithm for different values \(\alpha_k\) for the neo-Hookean-based damage model function \(\tilde{\Phi}\).
	}
	\label{fig:Erros_NH}
\end{figure}

\begin{figure}[htbp]
	\centering	
	\input{Pictures/RealDamage2D/NH/Phi_NH_2d_ReLU_ak_9_9}
	\caption{
    Comparison of the polyconvex envelopes \(\Phi^{\pc}_\pred\) (averaged over ten network realisations) and \(\Phi^{\pc}_{\delta}\) for the neo-Hookean model with $\hat{\nu}_k=[9,9]$ and for \(\alpha_k \in \{0.25, 0.75\}\). 
    The plot shows the cross-section along the diagonal axis.  
    The standard deviation of the ten predictions is marked by $\sigma$.
    }
\label{fig:NH_curves}
\end{figure}

\Cref{fig:Erros_NH} illustrates the relative errors between the predictions and the polyconvex envelopes computed by the SVPC~LP algorithm for different values of $\alpha_k$, including values that are outside the training domain. 
These errors are computed on a uniform $100 \times 100$ discretisation of the domain $[\nu_{\min},\nu_{\max}]^2$.
Across all considered parameter sets, the errors remain consistently low, ranging between 1\% and 2\%. 
In particular, it is important to note that the hypothesis stating that the energies become independent of $\alpha_k$ for $\alpha_k\geq \alpha_{\infty}$ is verified and validates our choice of $\alpha_\infty=4$. 
Additionally, \cref{fig:NH_curves} provides an example of re-shifted predictions, i.e.~application of \eqref{eq:splitpc}. 
It has to be noted that the predictions are in good agreement with the polyconvex envelopes computed by the SVPC~LP algorithm from \cref{sec:SVPC}. 

For comparison, the computation of one polyconvex envelope on a \(100 \times 100\)-grid on the box \([0.55, 18]^2\) using the SVPC~LP algorithm of \cref{sec:SVPC}, with the same discretisation parameters as employed for the learning data generation, takes \(1\) hour on a single CPU (exploiting the symmetry due to \(\Pid\)-invariance) while its prediction via the neural network takes only \(0.05\) seconds, highlighting the benefits of using a neural-network compression approach for the polyconvexification in engineering applications.
In addition, the neural network implemented in this example contains \num{14800} parameters, comprising both weights and biases. 
As before, this number of parameters remains smaller than the storage required for only two polyconvex envelopes on a \(100 \times 100\) grid.
These aspects emphasise once again the benefits of using a neural-network-based representation of the polyconvex envelope for parameter-dependent families of functions for applications in engineering problems.

\section{Conclusion}
We have demonstrated the effectiveness of a neural network design in predicting polyconvex envelopes with high accuracy and computational efficiency.
Our results show that such neural networks can generalise well beyond the learning dataset, enabling fast, multi-query evaluations and real-time computations. 
Moreover, we have introduced a splitting strategy which decouples the isotropic damage problem from the previous time step state, thereby improving the feasibility and robustness of the training process. 
Future research will focus on extending this framework to predict not only the polyconvex envelopes but also their derivatives as well as the incorporation of determinant constraints into the neural networks, as relevant for engineering applications in computational mechanics.
The results presented in this paper pave the way for complex material simulations in real engineering contexts.

\section*{Acknowledgement}
Many of the ideas in this paper were initially formulated and tested during Helena Althoff's Master's thesis \cite{Althoff2024}. Fruitful discussions with Daniel Balzani and Maximilian Köhler are gratefully acknowledged, and some ideas are the result of joint discussions with David Wiedemann. 

\printbibliography


\end{document}